\documentclass[3p]{elsarticle} 
\pdfoutput=1

\makeatletter
\def\ps@pprintTitle{%
 \let\@oddhead\@empty
 \let\@evenhead\@empty
 \def\@oddfoot{\centerline{\thepage}}%
 \let\@evenfoot\@oddfoot}
\makeatother

\usepackage{geometry}
\usepackage{setspace}

\usepackage[margin=1cm]{caption}
\usepackage{bbm}
\usepackage{soul}
\usepackage{eurosym}                           
\usepackage[latin1]{inputenc}                  
\usepackage{url}                               
\usepackage{longtable}                         
\usepackage{array}                             
\usepackage{graphicx,color}
\usepackage{amsthm}
\usepackage{amsbsy}

\usepackage{amssymb}
\usepackage{amsmath}
\usepackage{enumerate}
\usepackage{subfig}
\usepackage{graphicx,color}
\usepackage{epsfig}
\usepackage[ruled,vlined]{algorithm2e}
\usepackage{multirow,bigdelim}
\usepackage{multicol}
\usepackage{soul}

\usepackage{natbib}
\usepackage{scalerel}					 
\usepackage[english]{babel}
\newtheorem{theorem}{Theorem}
\newtheorem{definition}{Definition}
\newtheorem{remark}{Remark}

\graphicspath{{figures/}}

\theoremstyle{plain}

\theoremstyle{remark}

\theoremstyle{plain}

\theoremstyle{definition}

\newcommand{\e}{{\rm e}}        
\def\R{\mathbb{ R}}             
\def\E{\mathbb{ E}}             
\def\P{\mathbb{ P}}             
\def\F{\mathcal{F}}             
\renewcommand{\d}{{\rm d}}      
\def\dW{{\rm d}W}               
\def\dY{{\rm d}Y}               
\def\dt{{\rm d}t}

\def\dx{{\rm d}x}

\def\1{{\mathbbm{1}}}            

\def\texpsfig#1#2#3{\vbox{\kern #3\hbox{\includegraphics{#1}\kern #2}}\typeout{(#1)}}
\theoremstyle{plain}

\date{}

\newcommand{\zbox}[1]{
\noindent
\begin{center}
\framebox[14.5cm]{
\begin{minipage}{14cm}
#1
\end{minipage}
}
\end{center}
 }

\numberwithin{equation}{section}	     

\title{The Seven-League Scheme: Deep learning for large time step Monte Carlo simulations of stochastic differential equations}

\begin{document}

\author{Shuaiqiang Liu$^{1}$, Lech A. Grzelak$^{1,3}$, Cornelis W. Oosterlee$^{1, 2}$}


\address{$^1$Applied Mathematics (DIAM), Delft University of Technology, Delft, the Netherlands \\
$^2$Centrum Wiskunde $\&$ Informatica (CWI), Amsterdam, the Netherlands \\
$^3$Rabobank, Utrecht, the Netherlands}

\begin{abstract}
We propose an accurate data-driven numerical scheme to solve Stochastic Differential Equations (SDEs), by taking large time steps.
The SDE discretization is built up by means of  a polynomial chaos expansion method, on the basis of accurately determined stochastic collocation (SC) points.
By employing an artificial neural network to learn these SC points, we can perform Monte Carlo simulations with large time steps.  Error analysis confirms that this data-driven scheme results in accurate SDE solutions in the sense of strong convergence, provided the learning methodology is robust and accurate. With a method variant called the compression-decompression collocation and interpolation technique, we can drastically reduce the number of neural network functions that have to be learned, so that computational speed is enhanced.
Numerical experiments confirm a high-quality strong convergence error when using large time steps, and the novel scheme outperforms some classical numerical SDE discretizations. Some applications, here in financial option valuation, are also presented.
\end{abstract}


\maketitle

\section{Introduction}  \label{sec:introduction}

{\let\thefootnote\relax\footnotetext{The views expressed in this paper are the personal views of the authors and do not necessarily reflect the views or policies of their current or past employers.}}

The highly successful deep learning paradigm~\citep{2015nature} receives a lot of attention in science and engineering, in many different forms and flavors. 
Within numerical mathematics, the machine learning methodology has successfully entered the field of numerically solving   partial differential equations (PDEs)~\citep{data-driven-pde2019,2018SDE-DNN,Han8505,PiNN2019,SIRIGNANO20181339,TempoGAN2018}.
The aim with machine learning is to either speed up the solution process or  to solve high-dimensional problems that are not easily handled by the traditional numerical methods.
 There are essentially two types of deep learning approaches,  to approximate the solution of PDEs (see~\citep{2018SDE-DNN,Han8505,PiNN2019,SIRIGNANO20181339,TempoGAN2018}), and to design advanced numerical schemes, for example, the authors of \cite{data-driven-pde2019}, derived a data-driven high-order discretization to solve PDEs accurately on a coarse grid.

In this paper we will develop a highly accurate numerical discretization scheme for stochastic differential equations (SDEs), which is based on taking
possibly large discrete time steps. We ``learn'' to take large time steps, with the help of the Stochastic Collocation Monte Carlo sampler (SCMC) proposed by~\cite{scmc2019},  and by using an artificial neural network (ANN), within the classical supervised learning context.

SDEs are widely used to describe uncertain phenomena, in physics, finance, epidemics, amongst others, as a means to model and quantify uncertainty.
The corresponding solutions are stochastic processes.
Numerical approximation of the solution to an SDE is standard practice, as an analytic solution is typically not available.
The most commonly known technique to solve SDEs is based on Monte Carlo (MC) simulation, for which the SDE first needs to be discretized. 
There are quite a few applications, that could benefit from an accurate and efficient numerical method on the basis of a large time step discretization~\citep{2021largetime-step}, like, in finance, the valuation of path-dependent financial derivatives or financial risk management where counterparty credit risk play a role.

Basically, there are two ways to measure the convergence rate of discrete solutions to SDEs, by means of the approximation to the sample path or by approximation to the corresponding distribution.
This way, strong and weak convergence of a numerical SDE solution have respectively been defined (see~\cite{intro_sde1999}).
Weak convergence, the convergence in distributional sense, is often addressed in the literature. Moment-matching, for example, is a basic technique to improve weak convergence.
Strong, path-wise, convergence is particularly challenging, and requires accurate {\it conditional distributions}. There are natural approaches to improve strong convergence properties, i.e.  by adding higher order terms or by using finer time grids. However, these are nontrivial and costly, especially when considering multi-dimensional SDEs.

We aim to develop highly accurate numerical schemes by means of deep learning, for which the strong error of the discretization does not depend on the size of the simulation time step.
For this, we will employ the  SCMC method as
an efficient approach for approximating (conditional) distribution functions.
The distribution function of interest is then expanded as a polynomial in terms of a random variable which is {\it cheap to sample from} at given collocation points, and interpolation takes place between these points.
The resulting big time steps discretization, in which the SCMC methodology is combined with deep learning, is called {\em the Seven-League scheme}\footnote{With seven-league boots, we are marching through the time-wise direction, see also  \url{https://en.wikipedia.org/wiki/Seven-league_boot}
} here, and we abbreviate it by {\em the 7L scheme}.

There are different reasons to learn stochastic collocation points instead of the sample paths directly.  Stochastic collocation points have a specific physical meaning, which makes the data-driven scheme explainable.  Monte Carlo sample paths are random, while collocation points are path independent and deterministic (i.e. representing key features of a probability distribution), which simplifies the learning process  when using neural networks. Unlike the SCMC method, which  provides accurate Monte Carlo samples  given  a constant
time step and for one specific instance of the SDE parameters,  the 7L methodology enables us to generate samples for a wide range of time steps and for many different
instances of the model parameters (i.e. for a family of SDEs), by means of a neural network to learn the evolution of the collocation points over time for many different model parameters.

The 7L scheme is composed of two separate phases under the framework of supervised learning, i.e., an off-line (training) phase  and an on-line (prediction) phase. The training phase, which usually requires heavy computation and many data sets,  is done only once and offline. The prediction phase, which is a computationally cheap and highly efficient process, can be performed in an on-line fashion. 

The remainder of this paper is organized as follows. In Section~\ref{sec:sde}, SDEs, their discretization, stochastic collocation and the connection between SDE discretizations and the SCMC method are introduced.
In Section~\ref{sec:method}, the data-driven methodology is explained to address large time step simulation, i.e. the {\em 7L scheme}, for SDEs.
ANNs will be used  as function approximators to learn the stochastic (conditional) collocation points.  A brief description of their details is placed in Section~\ref{sec:ann}.
In Section~\ref{sec:cdc},  we introduce a decompression-compression technique to accelerate the computation.
This latter efficient variant is named the {\em 7L-CDC scheme} (i.e., seven-league compression-decompression scheme).
Section~\ref{sec:results} presents numerical experiments to show the performance of the proposed approach.
Furthermore, the corresponding error is analyzed.
Section~\ref{cenou} concludes.
\section{Stochastic differential equations and stochastic collocation}\label{sec:sde}

We first describe the basic, well-known SDE setting, and explain our notation.
\subsection{SDE basics}\label{sc:sdeb}
We work with  a real-valued random variable $Y(t)$, defined on the probability space $(\Omega,\Sigma,\mathbb{P})$ with filtration $\mathcal{F}_{t\in[0,T]}$, sample space $\Omega$, $\sigma$-algebra $\Sigma$ and probability measure  $\mathbb{P}$.
For the time evolution of $Y(t)$, consider the generic scalar It\^o SDE,
\begin{equation}\label{1a}
	\dY(t) = a(t,Y(t),{\boldsymbol\theta})\dt + b(t,Y(t),{\boldsymbol\theta}) \dW(t),  \;\;\;\;\; 0\leq t \leq T,
\end{equation}
with the drift term $a(t,Y(t),\boldsymbol\theta)$, the diffusion term $b(t,Y(t),\boldsymbol\theta)$, model parameters $\boldsymbol\theta$, Wiener process $W(t)$, and given initial value $Y_0:=Y(t=0)$.
When the drift and diffusion terms satisfy some regularity conditions (e.g., the global Lipschitz continuity~\citep[p.289]{karatzas1988brownian}), existence and uniqueness of the solution of~(\ref{1a}) are guaranteed.  The cumulative distribution function of $Y(t)$, $t\in[0,T]$, $F_{Y(t)}(\cdot)$, is available and the corresponding density function, evolving over time, is described by the Fokker-Planck equation~\citep{risken1984fokker}.

With a discretization in time interval $[0,T]$, $t_i=i\cdot T/N, \; i=0,\ldots N,$ with equidistant time step $\Delta t = t_{i+1}-t_i$, the discrete random variable at time $t_i$ is denoted by $Y(t_i)$. Traditional numerical schemes have been designed based on It\^o's lemma, in a similar fashion as the  Taylor expansion is used to discretize deterministic ODEs and PDEs.
The basic discretization, for each Monte Carlo path, is the Euler-Maruyama scheme~\citep{intro_sde1999}, which reads,
\begin{equation}
    \hat{Y}_{i+1} = \hat{Y}_{i} + a(t_i, \hat{Y}_i,\boldsymbol\theta) \Delta t+ b(t_i, \hat{Y}_i,\boldsymbol\theta) \sqrt{\Delta t} \hat{X}_{i+1},
\label{eq:euler}
\end{equation}
where $\hat{Y}_{i+1}:=\hat{Y}(t_{i+1})$ is a realization (i.e., a number) from random variable  $\tilde{Y}(t_{i+1})$, which represents the numerical approximation to exact solution $Y(t_{i+1})$ at time point $t_{i+1}$, and a realization $\hat{X}_{i+1}$ is drawn from the random variable $X$, which here follows the standard normal distribution $\mathcal{N}(0,1)$. Moreover, $\check{Y}(t_i)$ (a number) will be used as the notation for a realization of $Y(t_i)$.

In addition, the Milstein discretization~\citep{Milstein1975} reads,
\begin{equation}
    \hat{Y}_{i+1} =  \hat{Y}_{i} + a(t_i, \hat{Y}_i,\boldsymbol\theta) \Delta t+ b(t_i, \hat{Y}_i,\boldsymbol\theta) \sqrt{\Delta t} \hat{X}_{i+1} + \frac{1}{2} b'(t_i, \hat{Y}_i, \boldsymbol\theta) b(t_i, \hat{Y}_i,\boldsymbol\theta)  \Delta t(\hat{X}_{i+1}^2 -1),
\label{eq:milstein}
\end{equation}
where $b'(t_i, \cdot,\boldsymbol\theta)$ represents the derivative with respect to $\hat{Y}$ of $b(\cdot,\hat{Y},\boldsymbol\theta)$.
When the drift and diffusion terms are independent of time $t$, the SDE is called time-invariant.

Two error convergence criteria are commonly used to measure the SDE discretization accuracy, that is, the
convergence in the weak and strong sense.
Strong convergence, which is of our interest here, is defined as follows.
\begin{definition}
Let ${Y}(t_i)$ be the exact solution of an SDE at time $t_i$, its discrete  approximation $\tilde{Y}(t_i)$  with time step $\Delta t \in \mathbb{R}^+$ converges in the strong sense, with order $\beta_s \in \mathbb{R}^+$, if there exists a constant $K$ such that
\begin{equation}
        \E|\tilde{Y}(t_i) - {Y}(t_i)| \leq K(\Delta t)^{\beta_s}.
	\label{sstr}
\end{equation}
\end{definition}

It is well-known that the Euler-Maruyama scheme~\eqref{eq:euler} has strong convergence $\beta_s=0.5$, while the Milstein scheme~\eqref{eq:milstein} has $\beta_s$ = 1.0.
When deriving high order schemes for SDEs, the rules of It\^o calculus must be respected~\citep{intro_sde1999}. As a result, there will be eight terms in a Taylor SDE scheme with $\beta_s$ = 1.5, and twelve with $\beta_s$ = 2.0, and the computational complexity increases. As a consequence, higher order schemes are involved and somewhat expensive.
Convergence of the numerical solution for $\Delta t \rightarrow 0$ is guaranteed, but the computational costs increase significantly to achieve accurate solutions.

The generic form of the above mentioned numerical schemes to solve the It\^o SDE is as follows,
\begin{equation} \label{eq:ddscheme}
    \hat{Y}_{i+1} | \hat{Y}_i=  \sum_{j=0}^{m-1}\alpha_{j} \hat{X}_{i+1}^{j},
\end{equation}
where $m$ represents the number of polynomial terms, the coefficients $\alpha_j$ are pre-defined and equation-dependent.
For example, for the Euler-Maruyama scheme~\eqref{eq:euler}, with $m=2$, we have
\begin{equation} \label{eq:euler_coef}
\begin{cases}
\alpha_0  = \hat{Y}_{i} + a(t_i,\hat{Y}_i,\boldsymbol\theta) \Delta t, \\
\alpha_1 =b (t_i, \hat{Y}_i,\boldsymbol\theta) \sqrt{\Delta t}, \\
\end{cases}
\end{equation}
while for the Milstein scheme, with $m=3$, it follows that
\begin{equation}\label{eq:milstein_coef}
\begin{cases}
\alpha_0  = \hat{Y}_i + a(t_i, \hat{Y}_i,\boldsymbol\theta) \Delta t  + \frac{1}{2}b'(t_i, \hat{Y}_i, \boldsymbol\theta) b(t_i,\hat{Y}_i,\boldsymbol\theta), \\
\alpha_1 =b (t_i, \hat{Y}_i,\boldsymbol\theta) \sqrt{\Delta t}, \\
\alpha_2 = \frac{1}{2} b'(t_i,\hat{Y}_i,\boldsymbol\theta) b(t_i, \hat{Y}_i,\boldsymbol\theta). \\
\end{cases}
\end{equation}
With these explicit coefficients we arrive at the probability distribution of the random variable,
\begin{equation} \label{eq:distr_Yt}
 Y(t_{i+1}) | Y(t_i) \approx \tilde{Y}(t_{i+1}) |\tilde{Y} (t_i) \stackrel{d}{=}   \sum_{j=0}^{m-1}\alpha_{j} X^{j}.
\end{equation}
These discrete SDE schemes are based on a series of transformations of the previous realization to approximate the conditional distribution,
\begin{equation} \label{eq:condistr}
	\P \big[Y(t+\Delta t)<y|Y(t) \big]=F_{Y(t + \Delta t)|Y(t)}(y) \approx  F_{\tilde{Y}(t+\Delta t)|\tilde{Y}(t)}(y).
\end{equation}
A numerical scheme is thus essentially based on conditional sampling of $Y(t+\Delta t)|Y(t)$.
The Euler-Maruyama scheme draws from a normal distribution, with a specific mean and variance, to approximate the distribution in the next time point, while the Milstein scheme combines a normal and a chi-squared distribution. Similarly, we can derive the stochastic collocation methods.

\subsection{ Stochastic collocation method}\label{sc:scmc}
Let's assume two random variables, $Y$ and $X$, where the latter one is cheaper to sample from (e.g.,  $X$ is a Gaussian random variable).
These two scalar random variables are connected, via,
\begin{equation}
    F_Y(Y) \stackrel{d}{=} U \stackrel{d}{=} F_X(X),
\end{equation}
where  $U\sim \mathcal{U}([0,1])$ is a uniformly distributed random variable, $F_Y(\bar y):=\P[Y \leq \bar y]$ and $F_X(\bar x):=\P[X \leq \bar x]$ are cumulative distribution functions (CDF).  Note that $F_X(X)$ and $F_Y(Y)$ are random variables following the same uniform distribution. $F_Y(\bar y_n)$ and $F_X(\bar x_n)$ are supposed to be strictly increasing functions, so that the following inversion holds true,
\begin{equation} \label{eq:xydist}
    \bar y_n = F^{-1}_Y(F_X(\bar x_n))=:g(\bar x_n)  .
\end{equation}
where $\bar y_n$ and $\bar x_n$ are samples (numbers) from $Y$ and $X$, respectively. The mapping function, $g(\cdot) =  F^{-1}_Y(F_X(\cdot))$, connects the two random variables and guarantees that $F_X(\bar x_n)$ equals $F_Y(g(\bar x_n))$, in distributional sense and also element-wise. The mapping function should be approximated,  i.e., $g(\bar x_n)\approx g_m(\bar x_n)$, by a function which is cheap. When function $g_m(\cdot)$ is available, we may generate ``expensive'' samples, $\bar y_n$ from $Y$, by using the cheaper random samples $\bar x_n$ from $X$.

The Stochastic Collocation Monte Carlo method (SCMC) developed in~\cite{scmc2019} aims to find an accurate mapping function $g(\cdot)$ in an efficient way.
The basic idea is to employ Equation~\eqref{eq:xydist} at specific collocation points and approximate the function $g(\cdot)$ by a suitable monotonic interpolation between these points.
This procedure, see Algorithm~I, reduces the number of expensive inversions $F^{-1}_Y(\cdot)$ to obtain many samples from $Y(\cdot)$.

\zbox{
{\bf Algorithm I: SCMC Method} \\
\label{alg1}
Taking an interpolation function of degree $m-1$ (with $m\geq2$, as we need at least two collocation points), as an example, the following steps need to be performed:
\begin{itemize}
	\item[1.]  Calculate CDF $F_X(x_j)$ on the points $(x_1,x_2, ..., x_m)$, that are obtained, for example, from Gauss-Hermite quadrature, giving $m$ pairs $({x}_j, F_X({x}_j))$;
	\item[2.]  Invert the target CDF   ${y}_j=F^{-1}_Y(F_X({x}_j))$, $j=1,\ldots, m$, and form $m$ pairs $({x}_j, {y}_j)$;
	\item[3.]  Define the interpolation function, $y=g_m(x)$, based on these $m$  point pairs $({x}_j,{y}_j)$;
    \item[4.] Obtain sample $\hat{Y}$ by applying the mapping function $\hat{Y}=g_m(\hat{X})$, where  sample $\hat{X}$ is drawn from $X$.
\end{itemize}
}

The SCMC method parameterizes the distribution function by imposing probability constraints at the given collocation points.  Taking the Lagrange interpolation as an example, we can expand function $g_m(\cdot)$ in the form of polynomial chaos,
\begin{equation} \label{eq:scmc-lang}
    Y \approx g_m(X) = \sum_{j=0}^{m-1}\hat\alpha_{j} X^{j} = \hat\alpha_0+\hat\alpha_1X  + ... + \hat\alpha_{m-1} X^{m-1}.
\end{equation}
Monotonicity of interpolation is an important requirement, particularly when dealing with peaked probability distributions.

The Cameron-Martin Theorem~\citep{Cameron1947AnnMath}  states that any distribution can be approximated by a polynomial chaos approximation based on the normal distribution, but also other random variables may be used for $X$ (see, for example,~\cite{scmc2019}).

Clearly, Equation~\eqref{eq:distr_Yt} can be compared to Equation~\eqref{eq:scmc-lang}, as a discretization scheme to approximate the realization in the next time point.

\section{Methodology} \label{sec:method}
For our purposes, given $Y(t)$, the conditional variable $Y(t+\Delta t)$ can be written as,
\begin{equation} \label{eq:condY}
	Y(t + \Delta t)|Y(t) \approx g_m(X) = \sum_{j=0}^{m-1}\hat\alpha_{j}X^{j},
\end{equation}
where the coefficients $\hat\alpha_j \equiv \hat\alpha_j(\hat{Y}_i, t_i, t_{i+1},\Delta t,\boldsymbol\theta)$ are now functions of realization $\hat{Y}_i$.
Equation~\eqref{eq:condY}, with large $m$-values, holds for any $\Delta t$, particularly also for large $\Delta t$. As such the scheme can be interpreted as an {\em almost exact simulation scheme} for an SDE under consideration. By the scheme in~\eqref{eq:condY} we can thus take large time steps in a highly accurate discretisation scheme.
More specifically, a sample from the known distribution $X$ can be mapped onto a corresponding unique sample of the conditional distribution $Y(t+\Delta t)$ by the coefficient functions.


There are essentially two possibilities for using an ANN in the framework of the stochastic collocation method, the first being to directly learn the (time-dependent) polynomial coefficients, $\hat\alpha_{j}$, in~\eqref{eq:condY}, the second to learn the collocation points, ${y}_j$. The two methods are equivalent mathematically, but the latter, our method of choice, appears more stable and flexible.
Here, we explain how to learn the collocation points, ${y}_i$, which is then followed by inferring the polynomial coefficients.
When the stochastic collocation points at time $t+\Delta t$ are known, the coefficients in~\eqref{eq:condY} can easily be computed.

An SDE solution is represented by its cumulative distribution at the collocation points, plus a suitable accurate interpolation $g_m(x)$.  In other words, the SCMC method forces the distribution functions (the target and the numerical approximation) to strictly match at the collocation points over time.
The collocation points are dynamic and evolve with time.

\subsection{Data-driven numerical schemes} \label{sec:CondSCMC}

Calculating the conditional distribution function  requires generating samples conditionally on previous realizations of the stochastic process.
Based on a general polynomial expression, the conditional sample, in discrete form, is defined as follows,
\begin{equation} \label{eq:condsam}
	\hat{Y}_{i+1}|\hat{Y}_i =  \sum_{j=0}^{m-1}\hat\alpha_{i+1,j} \left(\hat{Y}_{i},  t_i, t_{i+1}-t_i, \boldsymbol\theta \right) \hat{X}_{i+1}^{j},
\end{equation}
where $\Delta t = t_{i+1} -t_i$, and the coefficients  $\hat\alpha_{i+1,j}$ are functions of the variables $\hat{Y}_{i},  t_i, t_{i+1}-t_i, \boldsymbol\theta$, for example, see Formulas~\eqref{eq:euler_coef} and~\eqref{eq:milstein_coef}.

In the case of a Markov process, the future doesn't dependent on past values. Given $\hat Y(t_i)$, the random variable  $\hat{Y} (t_{i+1})$ only depends on the increment $Y (t_{i+1})-Y (t_i)$. The process has independent increments, and
the conditional distribution at time $t_{i+1}$ given information up to time $t_i$ only depends on the information at $t_i$.

Similar to these coefficient functions, the $m$ {\em conditional} stochastic collocation points at time $t_{i+1}$, $y_j(t_{i+1})|\hat{Y}_i$, with $j=0,\ldots,m-1$,  can be written as a functional relation,
\begin{equation}
  {y}_j (t_{i+1}) | \hat{Y}_i = H_j \left(\hat{Y}_i,  t_i, t_{i+1}-t_i, \boldsymbol\theta \right).
\label{eq:markov-Hfun}
\end{equation}

A closed-form expression for function $H(\cdot)$ is generally not available.
Finding the conditional collocation points can however be formulated as a regression problem.

It is well-known that neural networks can be utilized as universal function approximators~\citep{Cybenko1989}. We then generate random data points in the domain of interest and the ANN should ``learn the mapping function $H_j(\cdot)$'', in an off-line ANN training stage. The SCMC method is here used to compute the corresponding collocation points at each time point, which are then stored to train the ANN, in a supervised learning fashion (see, for example,~\cite{Goodfellow-et-al-2016}).

\subsection{The Seven-League scheme}
Next, we detail the generation of the stochastic collocation points to create the training data. Consider a stochastic process $Y(\tau)$, $\tau\in[0,\tau_{max}]$, where $\tau_{max}$ represents the maximum time horizon for the process that we wish to sample from. When the analytical solution of the SDE is not available (and we cannot use an exact simulation scheme with large time steps), a classical numerical scheme will be employed, based on {\em tiny constant time increments} $\Delta \tau=\tau_{i+1}-\tau_i$, a discretization in the time-wise direction with grid points $0 < \tau_1 <\tau_2 < \ldots <\tau_N \leq \tau_{max}$,
to generate a sufficient number of highly accurate samples at each time point $\tau_i$, to approximate the corresponding cumulative functions highly accurately.
With the obtained samples, we approximate the collocation points, as follows,
\begin{equation} \label{eq:enq1}
  \hat{y}_j(\tau_i)=F^{-1}_{\tilde{Y}(\tau_i)}(F_X({x}_j)) \approx F^{-1}_{Y(\tau_i)}(F_X({x}_j))
\end{equation}
where $\hat{y}_j(\cdot)$ represents the approximate collocation points of $Y(t)$ at time $\tau_i$, and ${x}_j$, $j=1,\dots,M_s$, are collocation points of variable $X$. For simplicity, consider $X\sim\mathcal{N}(0,1)$, so that the points ${x}_j$ are known analytically and do not depend on time point $\tau_i$.  In the case of a normal distribution, these points are known quadrature points, and tabulated, for example, in~\citep{scmc2019}. After this first step,  we have the set of collocation points, $\hat{y}_j(\tau_i)$, for $i=1,\dots,N$ and $j=1,\dots,M_s$.
Subsequently, the $\hat{y}_j(\cdot)$ from~\eqref{eq:enq1} are used as the ground-truth to train the ANN.

In the second step,  we determine the {\it conditional} collocation points.
For each time step $\tau_i$ and collocation point indexed by $j$, a {\em nested Monte Carlo simulation} is then performed to generate the conditional samples.
Similar to the first step,  we obtain the conditional collocation points from each of these sub-simulations using~\eqref{eq:enq1}. This yields the following set of $M_c$ conditional collocation points,
\begin{eqnarray}
\label{eqn2}
	\hat{y}_{k|j}(\tau_{i+1}):=\hat{y}_{k}(\tau_{i+1})|\hat{y}_j(\tau_i)=F^{-1}_{\hat{Y}(\tau_{i+1})|\hat{Y}(\tau_i)=\hat{y}_j(\tau_i)}\biggl(F_X(x_{k|j})\biggr),
\end{eqnarray}
where $x_{k|j}$ is a conditional collocation point, and  $i \in \{0, 1,\dots,N-1\}$, $j \in \{1,\dots, M_s\}$, $k \in \{1,\dots,M_c\}$.
Note that, in the case of Markov processes, the above generic procedure can be simplified by just varying the initial value $Y_0$ instead of running a nested Monte Carlo simulation. Specifically, we then set  $\hat{Y}_i=\hat{Y}_0$,  $\tau_i=\tau_0$  and $\tau_{i+1}=\tau_0+\Delta\tau$ to generate the corresponding conditional collocation points.

The inverse, $F^{-1}_{\hat{Y}(\tau_{i+1})|\hat{Y}(\tau_i)}(\cdot)$, is often not known analytically, and needs to be derived numerically.
An efficient procedure for this is presented in~\citep{LLech}.
Of course, it is well-known that
the computation of $F^{-1}(p)$ is equivalent with the computation of the quantile at level $p$.

We encounter essentially four types of stochastic collocation (SC) points: ${x}_j$ are called the original SC points, $\hat{x}_j$ are original conditional collocation points, $\tilde{y}_j$ are the marginal SC points, and $\hat{y}_k|\cdot$ are the conditional SC points. For example,  $\hat{y}_{k}|\hat{Y_i}$ is conditional on a realization $\hat{Y}_i$. When a previous realization happens to be a collocation point, e.g., $\hat{Y}_i=\hat{y}_j$, we have $\hat{y}_{k|j}:=\hat{y}_{k}|\hat{y}_j$.

\begin{figure}[htb!]
\centering
\subfloat[Sample paths by 7L]{\label{fig:7L-3D} \includegraphics[width=0.5\textwidth]{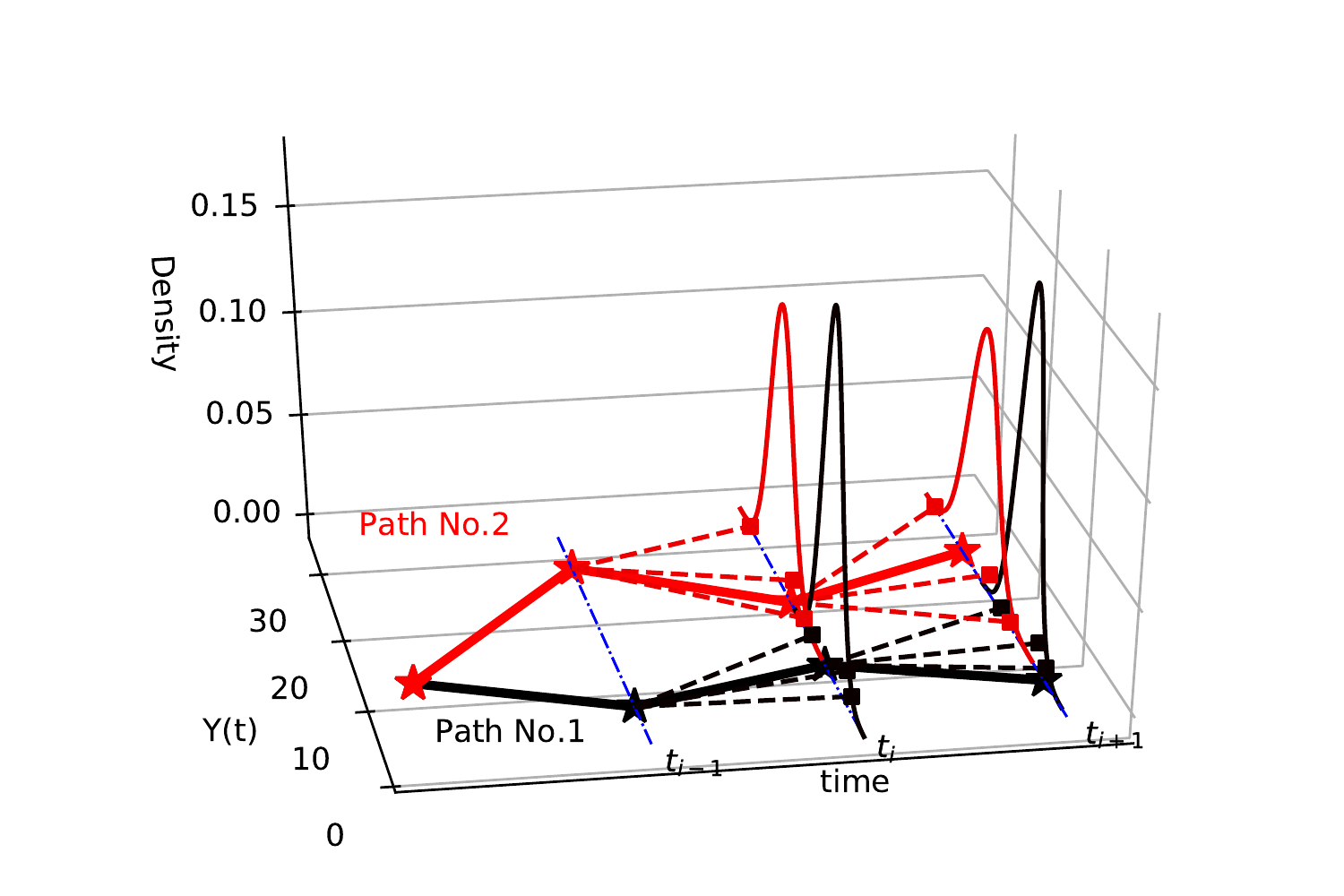}}
\subfloat[The 2D projection]{\label{fig:paths-7L-fig}{\includegraphics[width=0.5\textwidth]{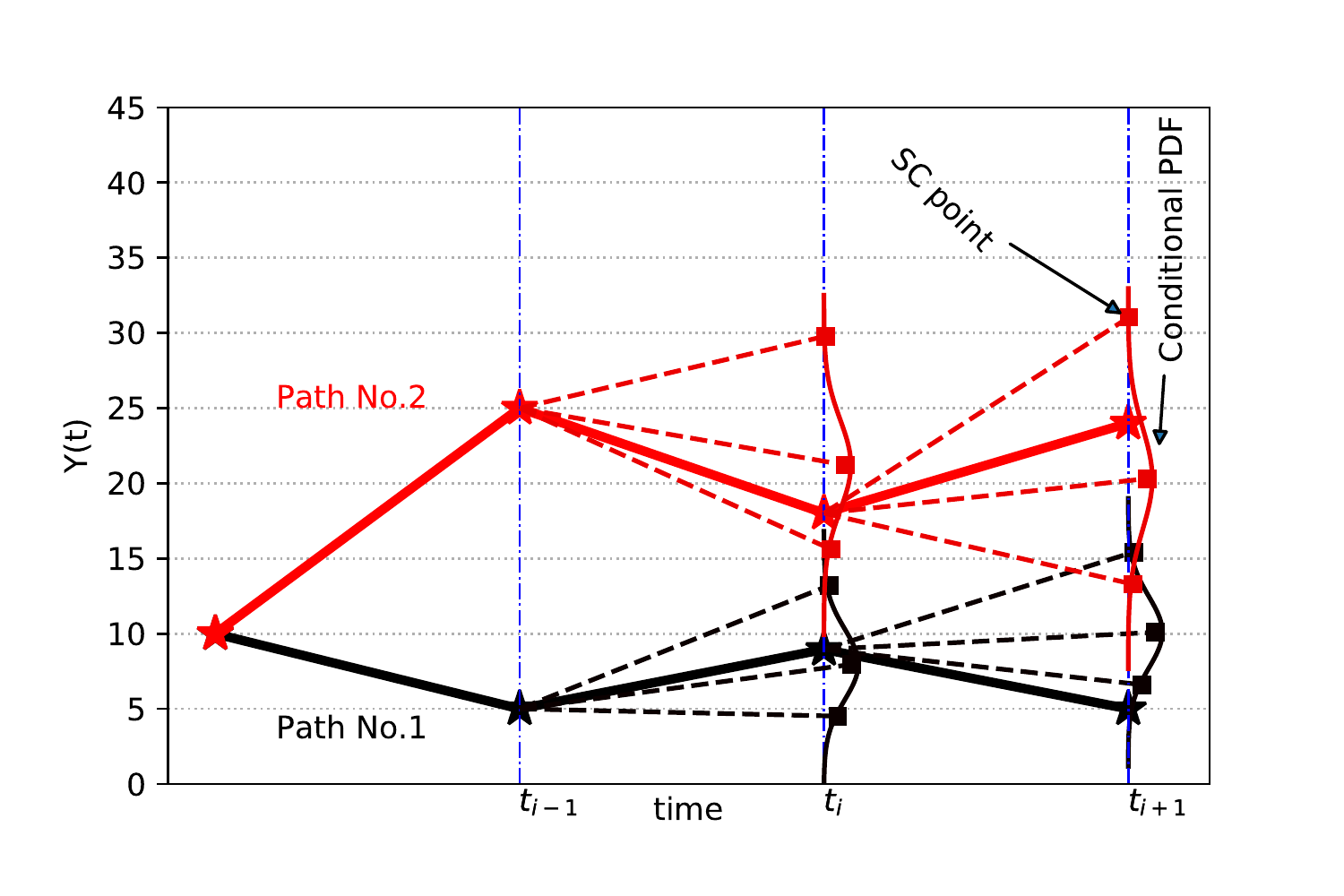}}}
\caption{ Schematic diagram of the 7L scheme. Left: Sample paths generated by 7L. Right: The 2D projection of Figure~\ref{fig:7L-3D}. Here conditional SC points, represented by $\blacksquare$, are conditional on a previous realization, denoted by $\bigstar$.  ``Conditional PDF'' is the conditional probability density function, defined by these conditional SC points. The density function, which is not required by 7L, is plotted only for illustration purposes.}
\label{fig:scheme-ann-scmc}
\end{figure}

When the data generation is completed, the ANNs are trained on the generated SC points to approximate the function $H$ in~\eqref{eq:markov-Hfun}, giving us a learned function $\hat{H}$. This is called the training phase.  With the trained ANNs, we can approximate new collocation points, and develop a numerical solver for SDEs, which is the Seven-League scheme (7L), see Algorithm~II. Figure~\ref{fig:scheme-ann-scmc} gives a schematic illustration of Monte Carlo sample paths that are generated by the 7L scheme.

\zbox{
	{\bf Algorithm II: 7L Scheme} \\
	\label{alg2}
\begin{enumerate}
      \item  Offline stage: Train the ANNs to learn the stochastic collocation points.
	      At this stage, we choose different ${\boldsymbol{\theta}}$ values, simulate corresponding Monte Carlo paths, with small constant time increments $\Delta \tau=\tau_{i+1}-\tau_i$ in $[0,\tau_{max}]$,
		generate the $\hat y_j$ and $\hat y_{k|j}$ collocation points, and learn the relation between input and output.
		So, we actually ``learn'' $H_k\approx \hat H_k$.
	      See Section~\ref{sec:ann} for the ANN details.

      \item  Online stage:
      Partition time interval $[0,T]$, $t_i=i\cdot T/N, \; i=0,\ldots N,$ with equidistant time step $\Delta t = t_{i+1}-t_i$.
      Given a sample $\hat{Y}_i$ at time $t_i$, compute $m$ collocation points at time $t_{i+1}$  using
      \begin{equation} \label{eq:funH}
        \hat{y}_j(t_{i+1}) | \hat{Y}_i = {\hat{H}_j} (\hat{Y}_i, t_i, t_{i+1}-t_i, \boldsymbol\theta ),  j=1,2,\ldots,m,
      \end{equation}
and form a vector ${\bf{\hat{y}}}_{i+1}=(\hat{y}_1(t_{i+1})| \hat{Y}_i,\hat{y}_2(t_{i+1})| \hat{Y}_i,\ldots,\hat{y}_m(t_{i+1})| \hat{Y}_i)$.
      \item Compute the interpolation function $g_m(\cdot)$, or calculate the coefficients  $\hat{\boldsymbol{\alpha}}_{i+1}$ (if necessary):
	      \begin{equation}
		      A({x}_{k|i+1}) \hat{\boldsymbol{\alpha}}_{i+1} = {\hat{\bf y}}_{i+1},
	    \label{eq:gm}
	      \end{equation}

		see Algorithm~I for details
		on the computation of original collocation points.
		We will compare monotonic spline, Chebyshev and the barycentric formulation of Lagrange interpolation for this purpose.
		See Section~\ref{sec:interp} for a detailed discussion.

\item Sample from $X$ and obtain a sample in the next time point, $\hat{Y}_{i+1}$, by $\hat{Y}_{i+1} | \hat{Y}_{i} = g_m(\hat{X}_{i+1})$,  or the coefficient form as follows,
       $$\hat{Y}_{i+1} | \hat{Y}_{i} =  \sum_{j=0}^{m-1}\hat{\alpha}_{i+1,j} \hat{X}_{i+1}^{j}.$$
\item Return to Step 2 by $t_{i+1} \xrightarrow{} t_{i}$, iterate until terminal time $T$.
\item Repeat this procedure for a number of Monte Carlo paths.
\end{enumerate}
}

\begin{remark}[Lagrange interpolation issue]
	In the case of classical Lagrange interpolation, matrix $A(x_{k|i})$ would be the Vandermonde matrix. In that case, it should not get too large, as the matrix would then suffer from ill-conditioning.
		However, when employing orthogonal polynomials, this drawback is removed.
		More details can be found in~\cite{scmc2019}.
\end{remark}

When the approximation errors  from ANN and SCMC are sufficiently small,
the strong convergence properties of the 7L scheme can be estimated, as follows,
\begin{equation}
	 \E|\tilde{Y}(t_i) - {Y}(t_i)| < \epsilon (\Delta \tau) \ll K(\Delta t)^{\beta_s},
\end{equation}
where  time step $\Delta \tau$ is used to define the ANN training data-set, and the actual time step $\Delta t$ is used for ANN prediction, with $\Delta \tau \ll \Delta t$.  Based on the trained 7L scheme, the strong error, $\epsilon(\Delta \tau)$, does thus not grow with the actual time step $\Delta t$. Particularly, let's assume $\Delta \tau= \Delta t/\kappa$, for example $\kappa=100$, when employing the Euler-Maruyama scheme with time step $\Delta \tau$ during the ANN learning phase, we expect a strong convergence of $O(\sqrt{\Delta \tau})$, which then equals $O(\sqrt{\Delta t/\kappa})$, while the use of the Milstein scheme during training would result in $O(\Delta t/\kappa)$ accuracy. When $\kappa=100$, the time step during the learning phase is 100 times smaller than $\Delta t$, which has a corresponding effect on the overall scheme's accuracy in terms of its strong, path-wise convergence.  {\em Moreover, the maximum value of the time step $\Delta t$ in the 7L scheme can be set up to $\tau_{max}$ for a Markov process. }  With a time step $\Delta t$,  we solve the SDE  in an iterative way until the actual terminal time $T$, which can be much larger than the training time horizon $\tau_{max}$. The detailed error analysis can be found in Section~\ref{sec:error-analysis}.

\subsection{The Artificial Neural Network} \label{sec:ann}

The ANN to learn the conditional collocation points is detailed in this subsection.
Neural networks can be utilized as powerful functions to approximate a nonlinear relationship.
In fact, we will employ a rather basic fully-connected neural network configuration for our learning task.


A fully connected neural network, without skip connections, can be described as a composition function, i.e.,
\begin{equation} \label{eq:fun-dnn}
\hat H(\Acute{\mathbf{x}}|\boldsymbol{\tilde{\Theta}}) = h^{(L)}(\ldots h^{(2)}(h^{(1)}(\Acute{\mathbf{x}};\boldsymbol{\tilde{\theta}}_{1});\boldsymbol{\tilde{\theta}}_{2});\ldots\boldsymbol{\tilde{\theta}}_{L_A}),
\end{equation}
where $\Acute{\mathbf{x}}$ represents the input variables, $\boldsymbol{\tilde{\Theta}}$ being the hidden parameters (i.e. weights and biases), $L_A$ the number of hidden layers.
We can expand the hidden parameters as,
 \begin{equation}
 \boldsymbol{\tilde{\Theta}} = (\boldsymbol{\tilde{\theta}}_1,\boldsymbol{\tilde{\theta}}_2,\ldots,\boldsymbol{\tilde{\theta}}_L) = (\mathbf{w}_1, \mathbf{b}_1, \mathbf{w}_2, \mathbf{b}_2,\ldots, \mathbf{w}_L,\mathbf{b}_{L_A}),
 \end{equation}
where $\mathbf{w}_\ell$ and $\mathbf{b}_\ell$ represent the weight matrix and  the bias vector, respectively, in the $\ell$-th hidden layer.

Each hidden-layer function,  $h^{(\ell)}(\cdot), \ell=1, 2, \ldots, L_A$, takes input signals from the output of a previous layer, computes an inner product of weights and inputs, and adds a bias. It sends the resulting value in an activation function to generate the output.

Let $z_{j}^{(\ell)}$ denote the output of the $j$-th neuron in the $\ell$-th layer. Then,
 \begin{equation}\label{eq:zz}
 z_{j}^{(\ell)}=\varphi^{(\ell)} \left(\sum_{i}w_{i,j}^{(\ell)}z^{(\ell-1)}_{i}+b_{j}^{(\ell)} \right),
\end{equation}
where $w_{i,j}^{(\ell)} \in \mathbf{w}_\ell$, $b_j^{(\ell)} \in \mathbf{b}_\ell$, and $\varphi^{(\ell)}$ is a nonlinear transfer function (i.e. activation function).  With a specific configuration,  including the architecture, the hidden parameters, activation functions and other specific operations (e.g., drop out), the ANN in~\eqref{eq:fun-dnn} becomes a deterministic, complicated, composite function.

Supervised machine learning~\citep{Goodfellow-et-al-2016} is used here to determine the weights and biases, where the ANN should learn the mapping from a given input to a given output, so that for a new input, the corresponding output will be accurately approximated. Such ANN methodology consists of basically two phases. During the (time-consuming, but off-line) training phase the ANN learns the mapping, with many in- and output samples,
while in the testing phase, the trained model is used to very rapidly approximate new output values for other parameter sets, in the on-line stage.

 In a supervised learning context, the loss function measures the distance between the target function and the function implied by the ANN. During the training phase, there are  many known data samples available,  which are represented by  input-output pairs $(\Acute{\mathbf{X}},\Acute{\mathbf{Y}})$.   With a user-defined loss function $\mathrm{L}(\boldsymbol{\tilde{\Theta}})$, training neural networks is formulated as
\begin{equation} \label{eq:argmin_dnn}
\arg \min_{\boldsymbol{\tilde{\Theta}}} \mathrm{L}(\boldsymbol{\tilde{\Theta}} | (\Acute{\mathbf{X}},\Acute{\mathbf{Y}})),
\end{equation}
where the hidden parameters are estimated to approximate the function of interest in a certain norm. More specifically, in our case, the input, $\Acute{\mathbf{x}}$, equals $\{\hat{Y}_i, t_i, t_{i+1}-t_i, \boldsymbol\theta\}$,
and the output, $\Acute{\mathbf{y}}$, represents the  collocation points ${\bf{\hat{y}}}_{i+1}$, as in Equation~\eqref{eq:funH}. In the domain of interest $\Acute{\Omega}$, we have a collection of data points $\{\Acute{\mathbf{x}}_k\}$, $k=1,\ldots,M_D$, and their corresponding collocation points $\{ \Acute{\mathbf{y}}_k\}$,  which form a vector of input-output pairs  $(\Acute{\mathbf{X}},\Acute{\mathbf{Y}})=\{(\Acute{\mathbf{x}}_k,\Acute{\mathbf{y}}_k)\}_{k=1,\ldots,M_D}$. For example, using the $L2$-norm, the discrete form of the loss function reads,
\begin{equation}
    \mathrm{L}(\boldsymbol{\tilde{\Theta}} | (\Acute{\mathbf{X}},\Acute{\mathbf{Y}}))  =  \sqrt{\frac{1}{M_D}  \sum_{k=1}^{M_D} \left(\Acute{\mathbf{y}}_k -\hat{H}(\Acute{\mathbf{x}}_k|\boldsymbol{\hat{\Theta}})\right)^2} .
\end{equation}
A popular approach for training ANNs is to optimize the hidden parameters via back-propagation, for instance, using stochastic gradient descent~\citep{Goodfellow-et-al-2016}.

\section{ An efficient large time step scheme: Compression-Decompression Variant} \label{sec:cdc}

The 7L scheme employs the ANNs to generate the conditional collocation points for  all samples of a previous time point, see Figure~\ref{fig:paths-7L-fig}.
The extensive use of ANNs in the methodology has an impact on the method's computational complexity.

In order to speed up the data-driven 7L scheme procedure, we introduce a compression-decompression (CDC) variant, {\em in the on-line validation phase.} Please note that the off-line learning phase is identical for both variants.
The so-called 7L-CDC scheme, to be developed in this section, only uses the ANNs to determine the conditional collocation points for the optimal collocation points of a previous time point.  All other samples will be computed by means of accurate interpolation. The computational complexity is reduced when the chosen interpolation is computationally cheaper than using ANNs.

By the compression-decompression procedure,  Monte Carlo sample paths based on SDEs can be recovered from a 3D matrix. We then employ the 7L scheme procedure only to compute the entries of the encoded  matrices $C_i$ at time point $t_i$, which leads to a reduction of the computational cost in many cases.

Next, we will explain the process of recovering the sample paths from a known matrix $C$ using the decompression method.

\subsection{CDC Variant}

With a time discretization  $\{t_0, t_1,t_2,\dots,t_N\}$, we define a 3D matrix $\hat{C} = \{\hat{C}_0,\hat{C}_1,\dots,\hat{C}_{N-1} \}$, which consists of  $N\times(M_s+1)\times(M_c+2)$ entries in total.  Recall that $M_s$ represents the number of collocation points and $M_c$ the number of conditional collocation points.
$M_s$ and $M_c$ may vary with time points $t_i$ (in case of an adaptive scheme, for example), but we use constant values for $M_s$ and $M_c$.
For each time point $t_i$, we construct a 2D matrix $\hat{C}_i$,
\begin{eqnarray}\label{eq:bigC}
\hat{C}_i = \left(
\begin{array}{cccccc}
	-& -& \hat{x}_1&\hat{x}_2 &\dots&\hat{x}_{M_c}\\
{x}_1& \tilde{y}_1(t_i)&\hat{y}_{1|1}(t_i)&\hat{y}_{2|1}(t_i)&\dots&\hat{y}_{M_c|1}(t_i)\\
{x}_2& \tilde{y}_2(t_i)&\hat{y}_{1|2}(t_i)&\hat{y}_{2|2}(t_i)&\dots&\hat{y}_{M_c|2}(t_i)\\
\vdots&\vdots&\vdots&\vdots&\vdots&\vdots\\
{x}_{M_s}&\tilde{y}_{M_s}(t_i)&\hat{y}_{1|M_s}(t_i)&\hat{y}_{2|M_s}(t_i)&\dots&\hat{y}_{M_c|M_s}(t_i)\\
\end{array}
\right)_{(M_s+1)\times(M_c+2)},
\end{eqnarray}
with $x_i$, $i=1,\ldots, M_s$, the original SC points, $\hat{x}_k$, $k=1,\ldots,M_c$, the $k$-th original conditional SC points, and the conditional SC points $\hat{y}_{k|j}(t_i) = \hat{y}_{k}(t_{i+1})|\hat{y}_{j}(t_i)$.  We thus represent the conditional SC points, $\hat{y}_{k}(t_{i+1})|\hat{y}_j(t_i)$, by matrix elements $c_{i,j,k}$.
The two empty cells in~\eqref{eq:bigC} are not addressed in the computation.
Moreover, at the last time point, $t_N$, $\hat{C}_{N}$ is not needed.
\begin{figure}[!h]
  \centering
	\subfloat[Marginal and Cond. SC points]{    \includegraphics[width=0.45\textwidth]{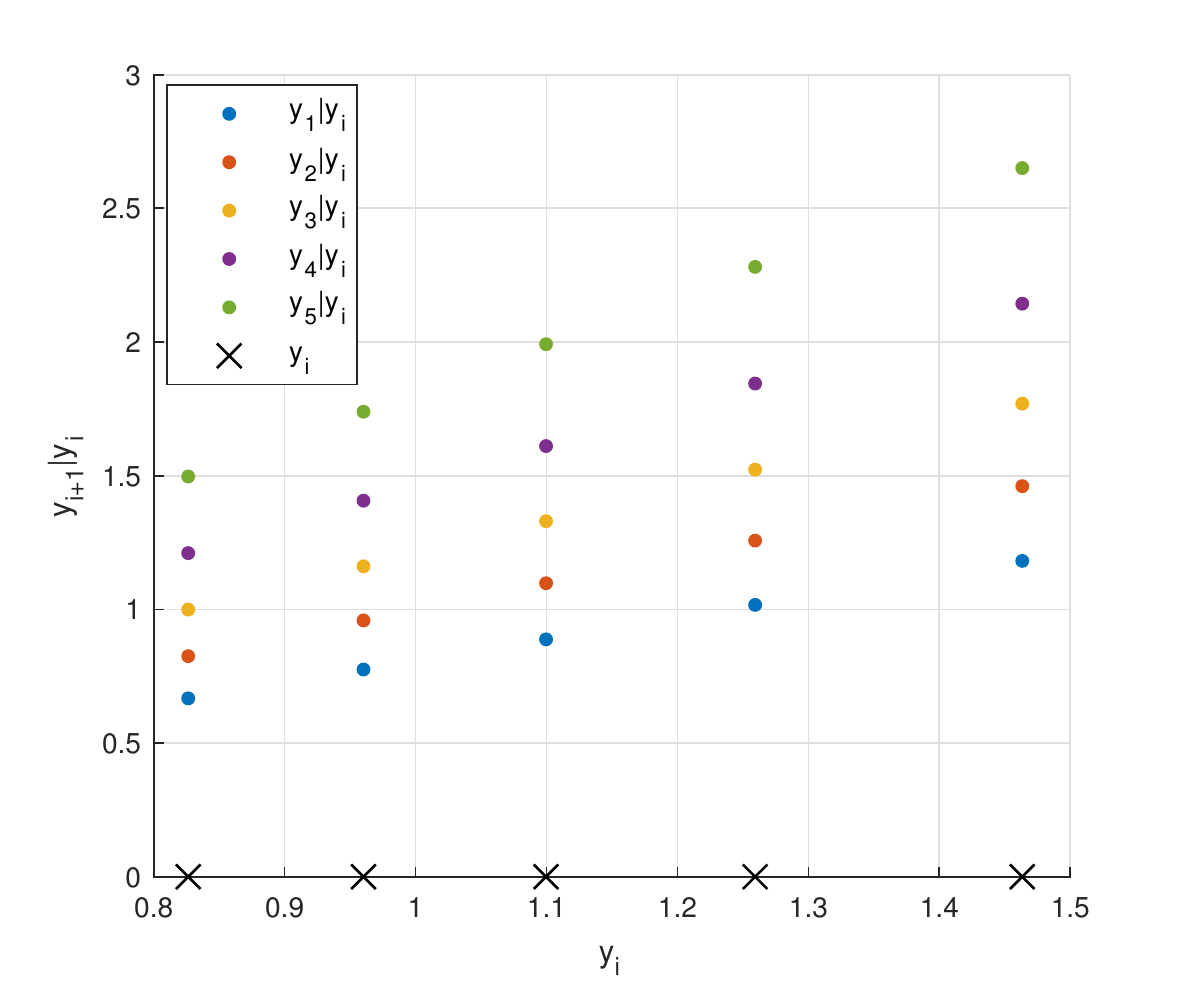}}
	\subfloat[Cond. SC points and densities]{\includegraphics[width=0.45\textwidth]{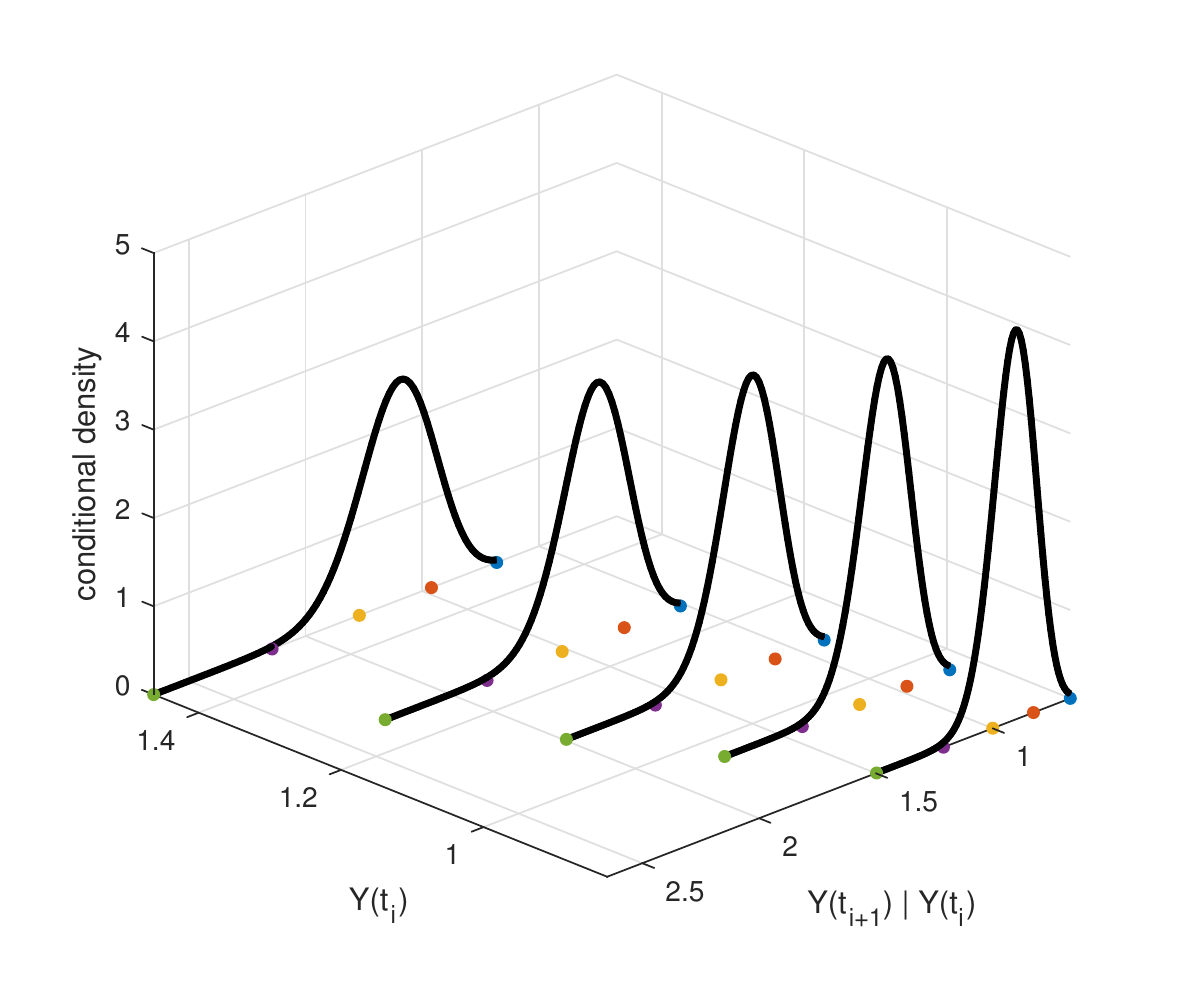}}
      \caption{Schematic illustration of matrix $C$, with five marginal SC points and five conditional SC points. The conditional SC points are dependent on the realization connected to the corresponding marginal SC point.  }
      \label{fig:collPoints3D}
\end{figure}

\begin{remark}[Time-dependent elements]
	As the original collocation points, ${x}_i$ and $\hat{x}_{k}$, do not depend on time, we can remove the first row and the first column of matrix $\hat{C}_i$ to obtain a time-dependent version, $C=\{C_0,C_1,\dots,C_{N-1}\}$, with the following elements,
\begin{eqnarray}\label{eq:smallC}
C_i=\left(
\begin{array}{ccccc}
 \tilde{y}_1(t_i)&\hat{y}_{1|1}(t_i)&\hat{y}_{2|1}(t_i)&\dots&\hat{y}_{M_c|1}(t_i)\\
 \tilde{y}_2(t_i)&\hat{y}_{1|2}(t_i)&\hat{y}_{2|2}(t_i)&\dots&\hat{y}_{M_c|2}(t_i)\\
\vdots&\vdots&\vdots&\vdots&\vdots\\
\tilde{y}_{M_s}(t_i)&\hat{y}_{1|M_s}(t_i)&\hat{y}_{2|M_s}(t_i)&\dots&\hat{y}_{M_c|M_s}(t_i)\\
\end{array}
\right)_{M_s\times (M_c+1)}.
\end{eqnarray}
\end{remark}

An entry in matrix $\hat{C}$ can be computed by the trained ANNs, as follows,
\begin{equation} \label{eq:cijk}
c_{i,j,k}: =\hat{y}_{k|j}(t_i)  = \hat{H}_k^{(M_c)}\Big(\hat{y}_{i,j},t_i, t_{i+1}-t_{i}, \boldsymbol\theta\Big),
\end{equation}
using the marginal SC points,
\begin{equation} \label{yij}
    \tilde{y}_j(t_i) =   \hat{H}_j^{(M_s)} \Big(Y_0,t_0, t_i-t_0, \boldsymbol\theta\Big),
\end{equation}
where $\hat{H}_j^{(\Lambda)} (\cdot), j=1,\dots,m$, $\Lambda=\{M_s,M_c\}$, represents the ANN function which approximates the $j$-th collocation point when $\Lambda=M_s$, and the $j$-th conditional collocation point when $\Lambda=M_c$.  When $M_c=M_s$, $\hat{H}_j^{(M_s)}=\hat{H}_j^{(M_c)}$. Figure~\ref{fig:collPoints3D} shows an example of the distribution of the conditional SC points when $M_c=3$ and $M_s=3$. When the matrices have been defined, all sample paths are compressed into a structured matrix. In other words, matrix $\hat{C}$ contains all the information needed to perform the Monte Carlo simulation of the SDEs, apart from the interpolation technique.

The resulting matrix $C$ will be {\em decompressed} to generate Monte Carlo sample paths with the help of an interpolation.
The process of decompression is straightforward given a matrix $\hat{C}$. In addition to the interpolation process $g_{M_c}(\cdot)$ in SCMC (see Equation~\eqref{eq:gm}), an interpolation $\tilde{g}(\cdot)$ is needed to compute conditional collocation points for previous realizations,  based on the matrix $\hat{C}$.

\begin{figure}[h!]
\subfloat[Paths for marginal SC points]{\includegraphics[width=0.45\textwidth]{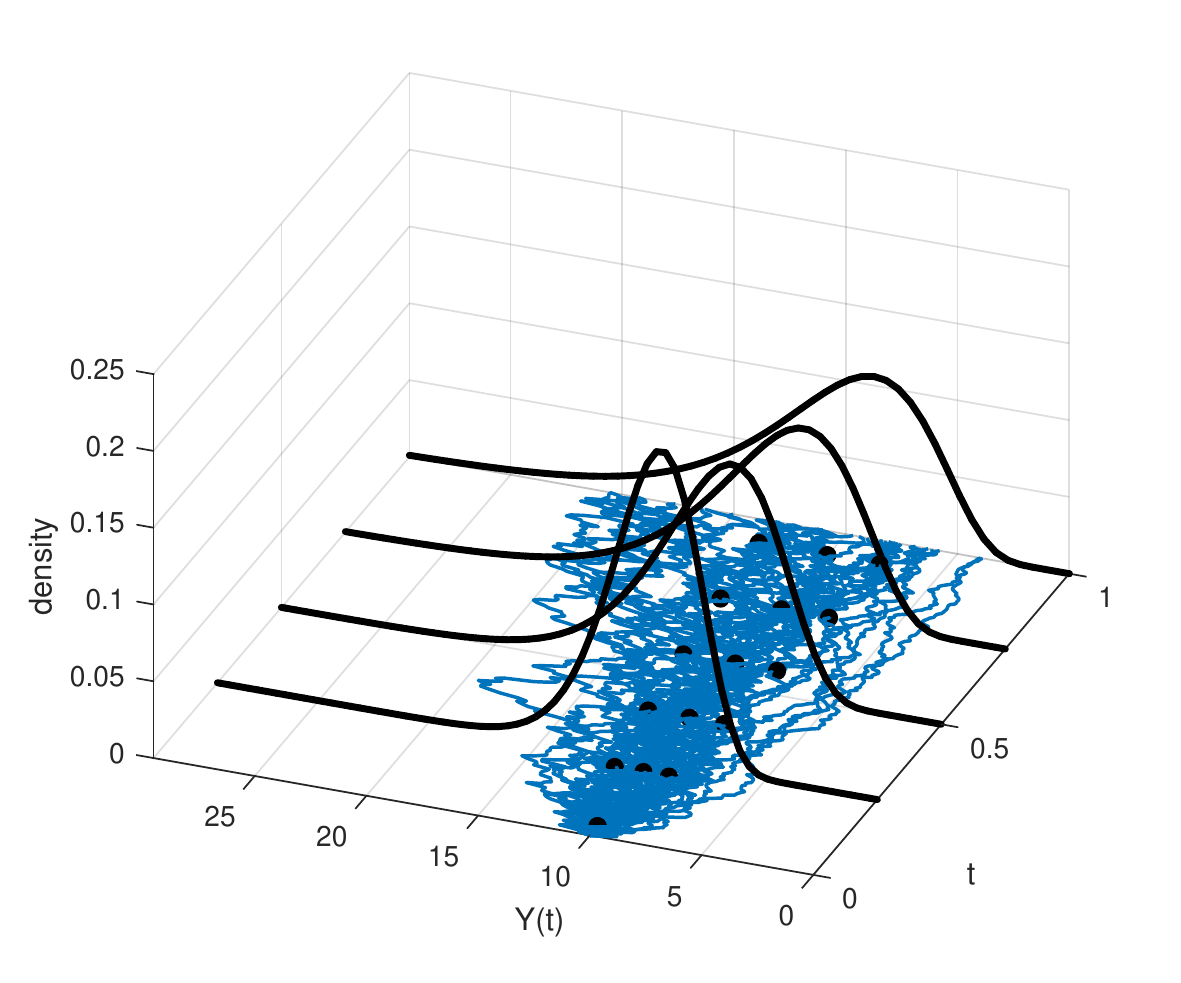}}
\subfloat[Sample paths by 7L-CDC]{\label{fig:paths-7L-cdc}{\includegraphics[width=0.55\textwidth]{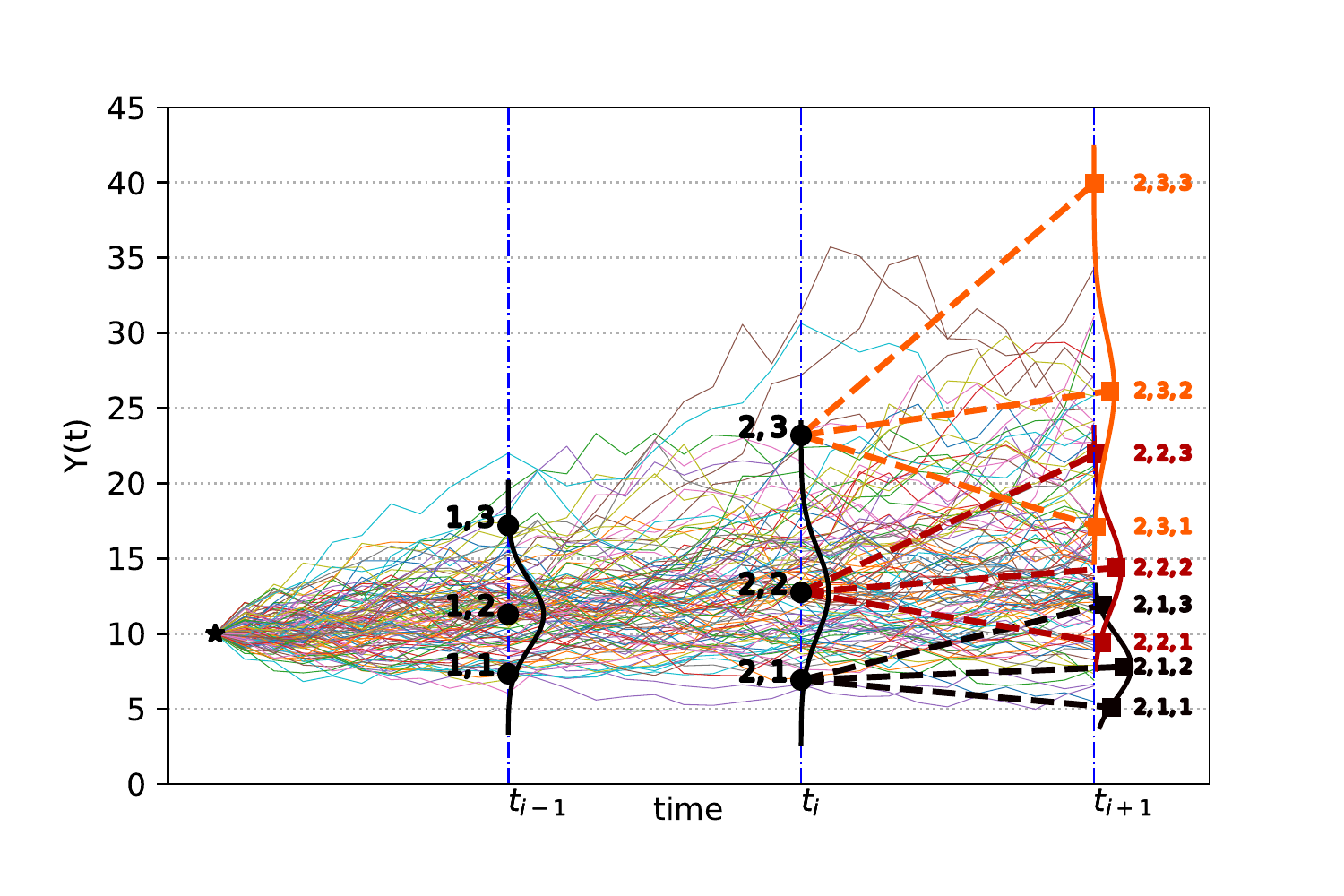}}}
\caption{ Schematic diagram of the 7L-CDC scheme at time $t_i$. Left: Marginal SC points, corresponding to Equation~\eqref{yij}. Right: Sample paths generated by 7L-CDC.  The triple \{2,1,3\}, in the picture, represents the third conditional SC point, dependent on the first marginal SC point at time point $t_2$. The above procedure is also applicable to other time points.}
\label{fig:scheme-ann-cdc}
\end{figure}

Suppose a vector of samples $\hat{\bf{Y}}_i$ at time $t_i$, and we wish to generate samples of $\hat{\bf{Y}}_{i+1}$.
For a specific sample $\hat{Y}^*_i$, we need to calculate $M_c$ conditional SC points. To obtain
the $k$-th ($1\leq k\leq M_c$) conditional SC point, we take marginal collocation points and their $k$-th conditional collocation points to form $M_s$ pairs $\{ ( \tilde{y}_1 (t_i),\hat{y}_{k|1} (t_i) ), (\tilde{y}_2 (t_i),\hat{y}_{k|2} (t_i) ),\dots,(\tilde{y}_{M_s},\hat{y}_{k|M_s} (t_i)) \}$.
This combination gives us the interpolation function $\hat{\bf y}=\hat g(\hat {\bf x})$.
Then  we can obtain the $k$-th conditional SC point of $\hat{Y}^*_i$,
\begin{equation} \label{eq:hatg}
    \hat{y}^*_k(t_{i+1})|\hat{Y}^*_i = \hat g(\hat{Y}^*_i).
\end{equation}

As a result, for each sample $\hat{Y}^*_i$, we obtain $M_c$ interpolation nodes, that form a set of pairs,  $\big(\hat{x}_1,\hat{y}^*_1(t_{i+1})|\hat{Y}^*_i \big)$,$\big(\hat{x}_2,\hat{y}^*_2(t_{i+1})|\hat{Y}^*_i \big)$, up to  $\big(\hat{x}_k,\hat{y}^*_{M_c}(t_{i+1})|\hat{Y}^*_i\big)$, which are used to determine the interpolation function $g_{M_c}(\cdot)$ required by SCMC.
Afterwards, to generate a new sample $\hat{Y}^*_{i+1}|\hat{Y}^*_i$,  the mapping function $g_{M_c}$ produces a conditional sample by taking in a random sample from $X$,
\begin{equation*}
    \hat{Y}^*_{i+1}|\hat{Y}^*_i= g_{M_c}(\hat{X}_{i+1}).
\end{equation*}

The choice of the appropriate number of (conditional) collocation points is a trade-off between the computational cost and the required accuracy. When the number of collocation points tends to infinity, the 7L-CDC scheme will resemble the 7L scheme from Section~\ref{sec:CondSCMC}.
A schematic picture is presented in Figure~\ref{fig:scheme-ann-cdc}.

\begin{remark}[Computation time]
During the on-line phase of the method, the total computation time of the large time step schemes consists of essentially two parts, calculation of the conditional SC points, and generating random samples by interpolation (the second part).
The difference between the 7L and 7L-CDC schemes is found in the computation of the conditional SC points, the generation of the samples is identical for both schemes.

In this first part, for the 7L-CDC scheme, the work consists of setting up matrix $C$ by the ANNs and computing the conditional SC points by the interpolation. In matrix $C$, there are $ M_s \times M_c \times N $ elements that are computed by the ANNs, where $N$ represents the number of time points, $M_s$ the number of collocation points and $M_c$ the number of conditional collocation points. Based on the $M_s$ collocation points, the interpolation is based on $M_c$ conditional collocation points for each path. For the 7L scheme, $M \times M_c \times N$ elements, where $M$ is the total number of paths, are computed by the ANNs. The time ratio between the 7L-CDC and  7L schemes is found as
\begin{equation} \label{eq:timeratio}
    \gamma = \frac{ t_I M + t_AM_s }{t_A M} = \frac{t_I}{t_A} +\frac{M_s}{M},
\end{equation}
with $t_A$ the computational time of the ANN (i.e., the function $\hat{H}(\cdot)$), $t_I$ for the interpolation (i.e., the function $\hat g(\cdot)$ in~\eqref{eq:hatg}), which is a polynomial function of $M_s$.   Given the fact that the number of sample paths is typically much larger than the number of SC points $M\gg M_s$, $$\gamma \approx \frac{t_I}{t_A}.$$  When the employed interpolation is computationally cheaper than the ANNs, $\gamma<1,$  so that the 7L-CDC scheme needs fewer computations than
the 7L scheme.
\end{remark}

\subsection{Interpolation techniques} \label{sec:interp}
To define the function $g_m(x)$ in~\eqref{eq:condY} or $\hat{g}(x)$ in~\eqref{eq:hatg}, we will compare three different interpolation techniques.

A bijective mapping function is obtained by the {\em monotonic} Piecewise Cubic Hermite Interpolating Polynomial (PCHIP)~\citep{PCHIP1980}. Assuming there are multiple data points, $(x_k,y_k)$,
using,
$$
h_k:=x_{k+1}-x_k, \; d_k:= \frac{y_{k+1}-y_k}{x_{k+1}-x_k},
$$
the derivatives $f'_k$ at the points $x_k$ are  computed as a weighted average,
$$\frac{\hat{w}_1+\hat{w}_2}{f'_k} = \frac{\hat{w}_1}{d_{k-1}} + \frac{\hat{w}_2}{d_k},
\;\; \mbox{ if } d_k \cdot d_{k-1} > 0,
$$
where  $\hat{w}_1:= 2h_k+h_{k-1}$ and $\hat{w}_2:=h_k +2h_{k-1}$. At each data point the first derivative is guaranteed to be continuous, and a cubic spline is used to interpolate between  the data points.
If $d_k \cdot d_{k-1} \leq 0$, then $f'_k=0$,
PCHIP requires more computations than a Lagrange interpolation, but it results in a monotonic function which is generally advantageous.

	The convergence of the stochastic collocation method is not really dependent on the monotonicity of the mapping function, so an interpolation based on {\em Lagrange polynomials} is possible in practice.
The barycentric version of Lagrange interpolation~\citep{Barycentric2004}, our second interpolation technique, provides a rapid and stable interpolation scheme, which
is applied when using Lagrange interpolation in our numerical experiments.
	With help of the basic Lagrange interpolation expressions, however, we can conveniently perform theoretical analysis.

The third technique is based on choosing the interpolation points carefully (e.g., as the Chebyshev zeros) to achieve a stable interpolation. The {\em Chebyshev  interpolation}~\citep{rivlin1990chebyshev} is of the form,
\begin{equation}
    g_m(x) =\sum_{j=0}^{m-1}\alpha_j p_{j}(x)= \alpha_0 + \alpha_1p_1(x) + ...+\alpha_{m-1}p_{m-1}(x),
\end{equation}
where $p_{m-1}(x)$ are interpolation basis functions, here Chebyshev orthogonal polynomials,  up to degree $m-1$. The Chebyshev nodes in the interval $[x_a,x_b]$ are computed as,
$$\tilde{x}_{k} = x_a + \frac{1+ \cos(\frac{\pi k}{m-1})}{2}(x_b-x_a), k=0,1,\ldots,m-1.$$
When the polynomial degree increases, the Chebyshev interpolation retains uniform convergence. In financial mathematics, Chebyshev interpolation has been successfully used, for example, to compute parametric option prices and implied volatility in~\citep{gass2018chebyshev,glau2019chebyshev,glau2019improved}.
When the interpolation points are not Chebyshev nodes (e.g., Gauss quadrature points), the Chebyshev coefficients can be estimated by means of a least squares regression, which is also called the Chebyshev fit.  In such case, the coefficients in~\eqref{eq:condsam} can be explicitly computed, in contrast to the barycentric Lagrange interpolation.  The selection of a suitable interpolation technique depends on various factors, for instance, speed,  monotonicity, availability of coefficients.
These three interpolation methods will be compared in the numerical section.

\subsection{Path-wise Sensitivity}
Often in computations with stochastic variables, we wish to determine the derivatives of the variables of interest, the so-called pathwise sensitivities. This is generally not a trivial exercise in a Monte Carlo setting, see, for example, the discussions in~\citep{fastgreeks,SmokingAF,finbook2019,JAIN201995}.
With our new large time step schemes, we determine the pathwise sensitivities of the computed stochastic variables in a natural way, based on the available information in the (conditional) SC points and the interpolation.
In this section, we derive the pathwise sensitivity of the state variable $Y(t)$ with respect to model parameters $\theta$.

The first derivative with respect to parameter $\theta$ of the conditional distribution in Equation~\eqref{eq:condY} reads,
\begin{equation} \label{eq:ywrttheta}
    \frac{\partial Y(t_{i+1})}{ \partial \theta} = \frac{\partial g(X)}{ \partial \theta} \approx \frac{\partial }{\partial \theta} \left(\sum_{j=1}^m \hat{y}_j(t) p _j(X) \right) = \frac{\partial }{\partial \theta} \left(\sum_{j=1}^m \hat{H}_j p_j(x) \right) = \sum_{j=1}^m \left( \frac{\partial \hat{H}_j}{\partial \theta} p_j(X) \right),
\end{equation}
where $p_j(X)$ are basis functions, which do not depend on the model parameters. For the derivative $\frac{\partial \hat{H}_j}{\partial \theta}$ in~\eqref{eq:ywrttheta} at time $t_i$, the expression of the ANN~\eqref{eq:fun-dnn}, given the specific activation function, is available. So, the function $\hat{H}$ is analytically differentiable. As a result, $\frac{\partial \hat{H}_j}{\partial \theta}$ can be easily computed, by means of automatic differentiation in the machine learning framework. Thus, we arrive at the sensitivity of a sample path with respect to model parameters, as follows,
\begin{equation} \label{eq:ywrtthetapath}
  \frac{\partial \hat{Y}_{i+1}}{ \partial \theta}  = \sum_{j=1}^{m} \frac{\partial \hat{H}_j}{\partial \theta} p_j(\hat{X}_{i+1}).
\end{equation}

\section{Numerical experiments} \label{sec:results}
In this section with numerical experiments we will give evidence of the high quality of our numerical SDE solver, by analyzing first in detail its components.
For this purpose, we mainly focus on the Geometric Brownian Motion SDE, which reads,
\begin{equation} \label{eq:gbmq}
    \d Y(t) = \mu Y(t) \dt+\sigma Y(t) \dW(t), \;\;\;\;\; 0 \leq t \leq T.
\end{equation}
The model parameters are the constant drift and volatility coefficients, i.e., $\boldsymbol\theta=\{\mu, \sigma\}$,  and the initial value is given by $Y_0$.
For~\eqref{eq:gbmq} a continuous-time analytic expression for the asset price at time $t$ is available, i.e.,
\begin{equation}
	Y(t)  = Y_0e^{(\mu-\frac{1}{2} \sigma^2) (t-t_0)+\sigma ( W(t)-W(t_0)) } \stackrel{d}{=} Y_0e^{(\mu-\frac{1}{2} \sigma^2) (t-t_0)+\sigma \sqrt{t-t_0}X},
	\label{exa}
\end{equation}
where  $X\sim \mathcal{N}(0,1)$, and $Y(t)$ is governed by the lognormal distribution. The derivative of the stock price with respect to volatility $\sigma$ is available in closed form, and reads,
\begin{equation} \label{eq:gbmsense}
\frac{\partial Y(t)}{\partial \sigma} \stackrel{d}{=} Y(t)(-\sigma (t-t_0) + \sqrt{t-t_0}X).
\end{equation}
This expression will be used as the reference value of the sensitivity obtained from the 7L discretization.

Furthermore, the Ornstein-Uhlenbeck process is explained and also analyzed, in Subsection~\ref{ouou}.
We will employ the large time step discretization, in which the conditional collocation points are computed by the trained ANN,
and compare the results of the novel scheme with those obtained by the Milstein SDE discretization.
\subsection{ANN Training Details} \label{sec:trainprocess}
GBM and the OU process are Markov processes, so the conditional distribution at time $t_{i+1}$ given information up to time $t_i$ only depends on the information at time $t_i$.
The ANN~\eqref{eq:markov-Hfun} will therefore be used for the conditional collocation stochastic points,
with $\boldsymbol\theta = \{\mu, \sigma\}$, for GBM, and $\boldsymbol\theta = \{\overline{Y},\sigma,\lambda\}$ for the OU process (as will be discussed in Subsection~\ref{ouou}).

Regarding the size of the compression-decompression matrix, the more conditional collocation points, the better the accuracy of the 7L-CDC method. A 5x5 matrix size (i.e., five marginal and five conditional SC points) is preferred, taking into account the computing effort and the accuracy. In~\cite{scmc2019} it has been discussed and shown that highly accurate approximations could already be obtained with a small number of collocation points.

As the first method component, we evaluate the quality of the ANN which defines the collocation points, for the GBM dynamics.
For this purpose, $M_L$ random points (i.e. sets of input parameters) are generated by using Latin Hyper-cube Sampling (LHS) in the domain of interest for the three parameters $(Y_0,\mu, \sigma)$, see Table~\ref{table:ann_train_data}. As the second step, for each point a Monte Carlo method is employed to simulate the discretized SDE based on the tiny time step $\Delta \tau$.
We use an Euler-Maruyama time discretization for this purpose,
with $N_\tau$ the number of time points and the time horizon $\tau_{max} = N_\tau \cdot \Delta \tau$. At each time step, $j=1,...,N_\tau$, the conditional distribution function
$F_{Y(t_j)|Y_0, \mu,\sigma}(\cdot)$ is computed, based on the many generated MC paths.  This way, the resulting collocation points for the ``big time step'', $\Delta t = j\cdot\Delta \tau$, are also obtained, to form the required training data set.

We set $\tau_{max}=1.6$, $N_\tau=160$, $M_L=500$.
The amount of training data used is given by  $M_{train} = M_L \cdot N_\tau=80,000$ samples in total, which are divided into an ANN training ($90\%$) and an ANN testing ($10\%$) set.  

\begin{table}[!h]
\begin{center}
 \caption{ Training data, $\Delta \tau = 0.01$. Here is an example for training on five SC points. }
 \begin{tabular}{ c|c | c | c }
  \hline
  ANN & Parameters       & Value range   & Method \\ \hline
  \multirow{3}{*}{input} & drift, $\mu$     & (0.0, 0.10] & LHS \\
  &volatility, $\sigma$    & [0.05, 0.60] & LHS \\
  &value, $Y_0$     & [0.10, 15.0] & LHS \\
  &time, $\tau_{max}$    & (0.0, 1.60] & Equidistant\\

 \hline
 \multirow{1}{*}{$\hat{H}_1(\cdot)$ output} & point, $\hat{y}_1$ & $(0.0, 25.65)$ & SCMC\\
\multirow{1}{*}{$\hat{H}_2(\cdot)$ output} & point, $\hat{y}_2$ & $(0.0,  25.98)$ & SCMC\\
\multirow{1}{*}{$\hat{H}_3(\cdot)$ output} & point, $\hat{y}_3$ & $(0.0, 27.84)$ & SCMC\\
\multirow{1}{*}{$\hat{H}_4(\cdot)$ output} & point, $\hat{y}_4$ & $(0.0,  54.67)$ & SCMC\\
\multirow{1}{*}{$\hat{H}_5(\cdot)$ output} & point, $\hat{y}_5$ & $(0.0, 154.35)$ & SCMC\\   \hline
 \end{tabular}
 \label{table:ann_train_data}
 \end{center}
 \end{table}
 
The ANN hyper-parameters have an impact on the errors from optimization related to training the ANN, as well as on the model performance.
The approximation capacity does not only depend on  the number of hidden parameters, but also on the network structure (i.e., on the width and depth of the network). In principle, deep neural networks have more powerful expressiveness than shallow neural networks.
The fully connected neural network employed will be composed of one input layer, one output layer and four hidden layers. Each hidden layer consists of 50 neurons, with Softplus, i.e. $\varphi(x)= \ln(1+e^x)$ as the activation function~\citep{activation2018}. Before training the ANN,  the hidden parameters are initialized via the Glorot technique~\citep{glorot2010}.  Training goes in batches. At each iteration, the stochastic gradient based optimizer, Adam~\citep{Adam2014}, randomly  selects a portion of the training samples according to the batch size, to calculate the gradient for updating the hidden parameters. In an epoch, all training samples have been processed by the optimizer.  The mean squared error, which measures the distance between the ground-truth and the model values in supervised learning,  is used to update the hidden parameters during training.
 The measure MAE (Mean Absolute Error), i.e.,  $\text{MAE} = \frac{1}{M_{train}} \sum_j |{y}_j-\hat{y}_j|,$ is also estimated, as the path-wise error of the 7L scheme is related to the maximum absolute difference in the approximated collocation points, $\hat{y}_j$, in Section~\ref{sec:pathwise error}, see the derivation in the next section.

\begin{figure}[!h]
\centering
\subfloat[$\hat{y}_2$]{\label{fig:test_r2_y2} {\includegraphics[width=0.45\textwidth]{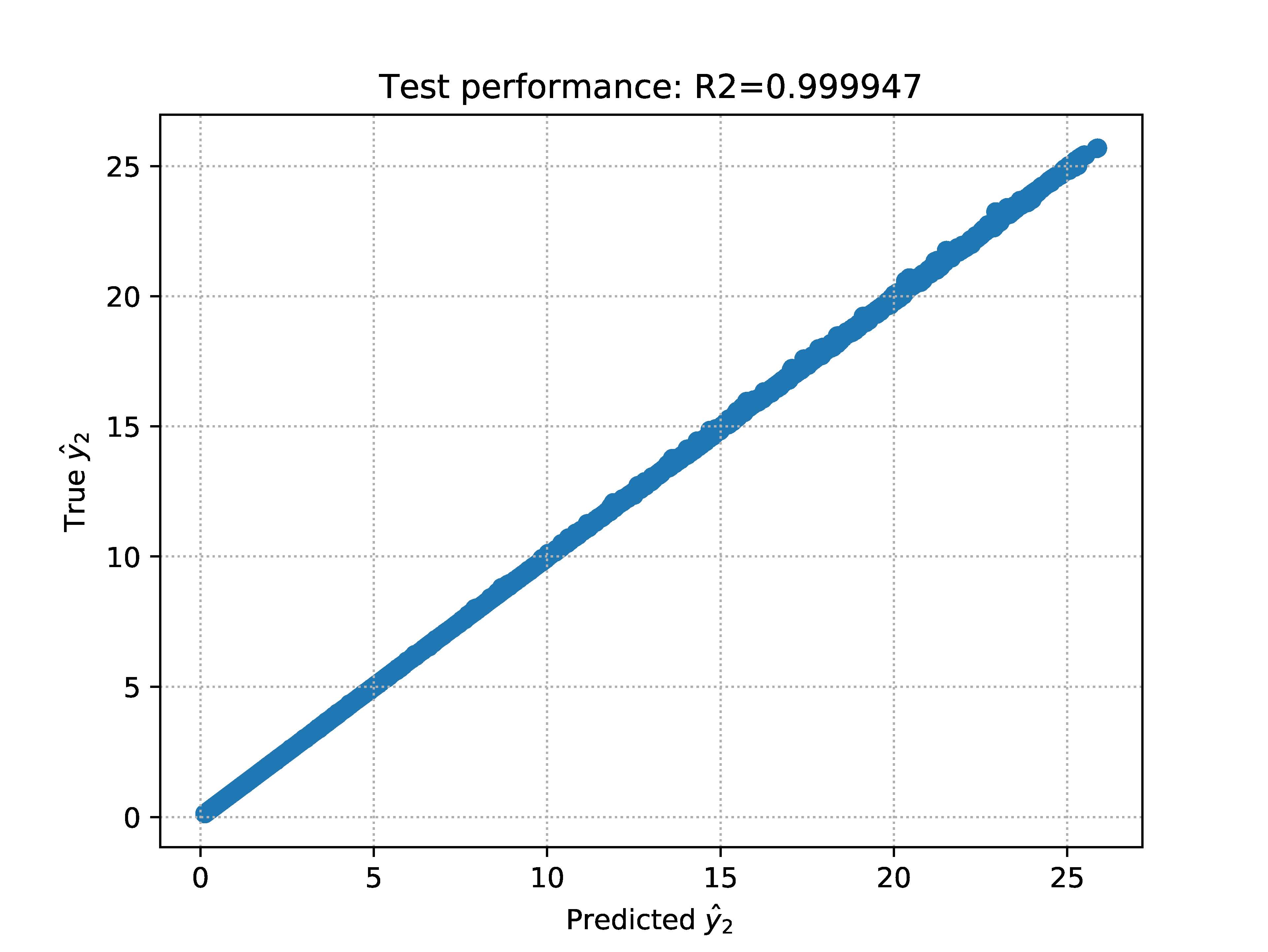}}}%
\subfloat[$\hat{y}_4$]{\label{fig:test_r2_y4}{\includegraphics[width=0.45\textwidth]{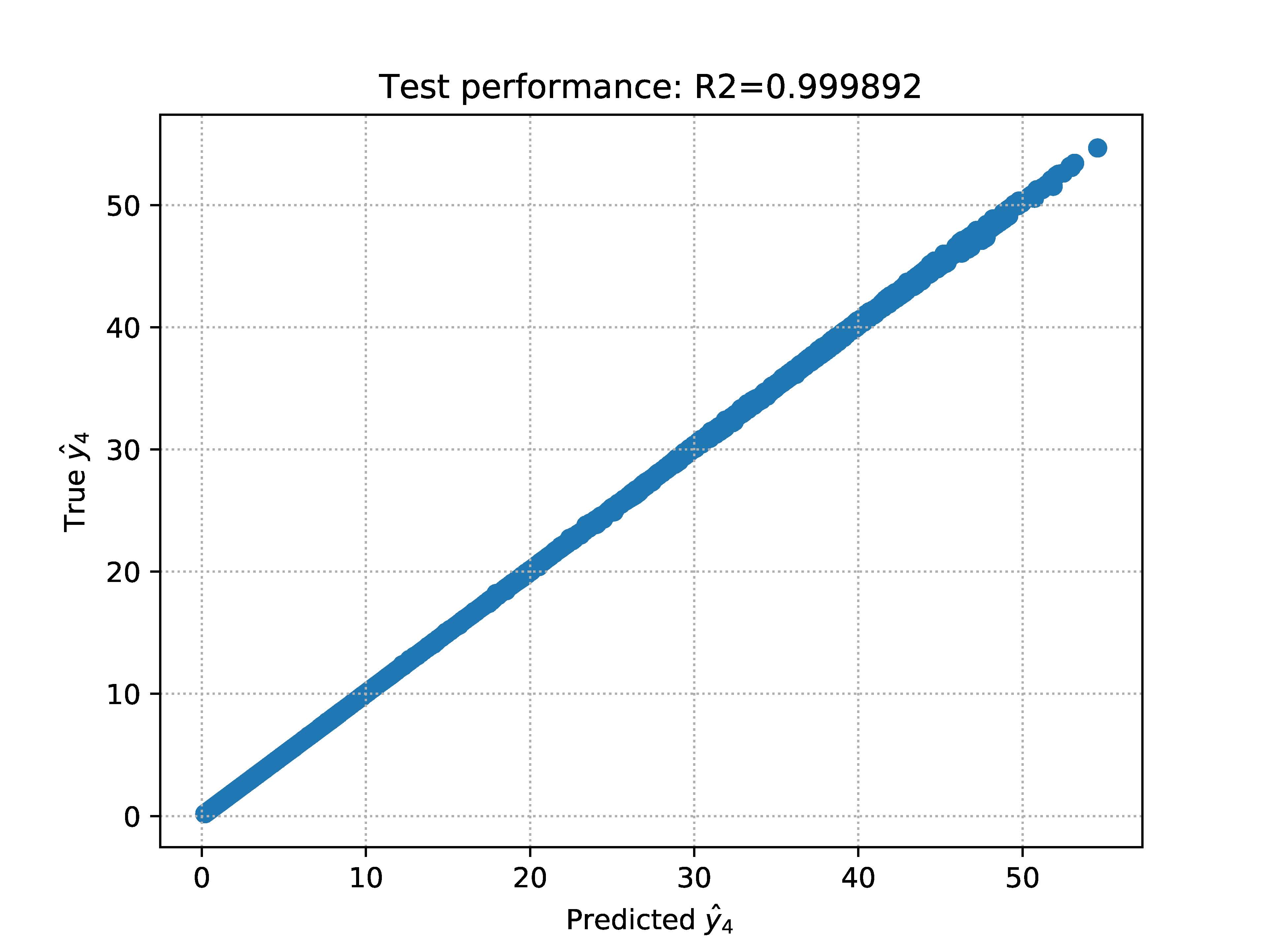}}}
\caption{The goodness of fit on test data set. Two scatter plots show the relation between the predicted values and the ground truth.}
\label{figure:r2-ANN-performance}
\end{figure}

The training process starts with a relatively large learning rate (i.e $10^{-3}$) to avoid getting stuck in local optima. After 1000 epochs, the learning rate is reduced to $10^{-4}$, followed by training 500 more epochs, to achieve a steady convergence. 
Afterwards, the trained ANN is evaluated on the testing data set, with the results presented in Figure~\ref{figure:r2-ANN-performance} (for two of the collocation points) and Table~\ref{table:measures_ann}. Clearly, the predicted values fit very well with the true values of the stochastic collocation points. This implies that the trained ANNs reach a highly satisfactory generalization, and generate accurate and robust approximation results for all five collocation points.

\begin{table}[htp]
\begin{center}
\caption{The approximation performance on test data set.}
\scalebox{1.0}{
   \begin{tabular}{  c | c | c | c | c | c}
    \hline
    SC points &  $\hat{y}_1$ & $\hat{y}_2$ & $\hat{y}_3$ & $\hat{y}_4$ & $\hat{y}_5$   \\ \hline
    $R^2$  & 0.999891 & 0.999947 & 0.999980 & 0.999892 & 0.999963 \\
    MAE & 0.026 &  0.027  & 0.021 & 0.071 & 0.066 \\
    \hline
  \end{tabular}  }
\label{table:measures_ann}
\end{center}
\end{table}

\subsection{Error analysis, the Lagrangian case} \label{sec:error-analysis}
There are essentially two approximation errors in the 7L scheme, a neural network approximation error when generating the collocation points, and an SCMC error when representing the conditional distribution function.

Considering $d$ inputs, the neural network may approximate any function $\zeta_{d,n}$, from the function space $C^{n-1}([0,1]^d)$, where the derivatives up to order $n-1$ are Lipschitz continuous~\cite{ErrorReLuDNN2017}. The input and output variables can be normalized to the unit interval $[0,1]$. With a fixed network architecture during training, the approximation error can be assessed, as follows.

\begin{theorem}
From~\cite{ErrorReLuDNN2017}. Given any $\hat{\epsilon} \in (0,1)$, there exists a neural network  which is capable of approximating any function $\zeta_{d,n}$ with error $\hat{\epsilon}$, based on the following configuration:
\begin{itemize}
    \item at least piece-wise activation functions,
    \item at least $\tilde{c}(\ln(1/\hat{\epsilon} )+1 )$ hidden layers and $\tilde{c}\hat{\epsilon}^{-d/n}(\ln(1/\hat{\epsilon}) +1 )$ weights and computation units, where $\tilde{c}:=\tilde{c}(d,n)$ depends on the parameters $d$ and $n$.
\end{itemize}
\end{theorem}
When the architecture is dynamic, the error bound can be further reduced, as shown in~\cite{NewErrorBoundsDNN2019} and~\cite{ErrorReLuDNN2017}.
One of the assumptions is that the ANNs are sufficiently trained, so that the optimization error is negligible.

The error from the SCMC methodology was derived in~\cite{scmc2019}. The optimal collocation points, ${x}_i$, $i=1,\dots,m$, correspond to the zeros of an orthogonal polynomial. In the case of Lagrange interpolation, when the collocation method can be connected to Gauss quadrature,
we have
\begin{eqnarray}
\int_\R \Psi(x)f_X(x) \dx=\sum_{i=1}^m \Psi({x}_i)
\omega_i+ \epsilon_m = \epsilon_m,
\end{eqnarray}
with $\Psi(x)=\left(g(x) - g_m(x)\right)^2$, the difference between the target and the SC approximated function, $f_X(x)$ the weight function, and $\omega_i$ the quadrature weights. When the Gauss-Hermite quadrature is used with $m$ collocation points.
the approximation error of the CDF can be estimated as,
\begin{eqnarray}\label{epsilon_scmc}
\epsilon_m=\frac{m!\sqrt{\pi}}{2^m}\frac{\Psi^{(2m)}(\xi_1)}{(2m)!},
\end{eqnarray}
where $ \xi_1 \in (-\infty, \infty)$ and the distance function $$\Psi(x)=\left(g(x) - g_m(x)\right)^2\approx \left(\frac{1}{m!}\left.\frac{\partial^m g(x)}{\partial x^m}\right|_{x=\xi_2}\prod_{k=1}^m (x-{x}_k)\right)^2,$$
with $\xi_2 \in [{x}_1,{x}_{m-1}]$.
In other words, the error of approximating the target CDF converges exponentially to zero when the number of corresponding collocation points increases.

At each time point $t_i$, the process $Y(t_i)$ is approximated  using the collocation method, by a polynomial $g_m(X)$, i.e., in the case of classical Lagrange interpolation,
using $\ell_j(\bar x) = p_j(\bar x)$,
\begin{equation}
	Y(t_i) \approx  \tilde{Y}(t_i)=g_m(X)=\sum_{j=1}^m {y}_j(t_i)\ell _j(X), \;\; \ell_j(\bar{x})=\prod_{k=1}^m\frac{X-{x}_k}{{x}_j-{x}_k},
\end{equation}
where the collocation points  ${y}_j(t_i)=F^{-1}_{Y(t_i)}(F_X({x}_j))$.  Because of the use of an ANN, the collocation points are not exact, but they are approximated with ${y}_j(t_i)-\hat{y}_j(t_i)=\epsilon_j^{A}$, where $\hat{y}_j(t_i)$ represents the ANN approximated value. The error associated with $\epsilon_j^{A}$ can be estimated as in~\citep{NewErrorBoundsDNN2019}. The impact of $\epsilon_j^A$ on the obtained output distribution needs to be assessed. Let $\tilde{g}_m$ denote the approximate function based on the predicted ANN collocation points $\hat{y}_j(t_i)$, and $x$ a random sample from the standard normal distribution $X$.  The approximation error, in the strong sense, is given by
\begin{eqnarray}
\E\left[ |g_m(x)-\tilde{g}_m(x) |\right]&=&\E\Big|\sum_{j=1}^m{y}_j(t_i)\ell _j(x)-\sum_{j=1}^m \hat{y}_j (t_i)\ell _j(x)\Big|\nonumber \\
&=&\int_\R\Big|\sum_{j=1}^m {y}_j(t_i)\ell _j(x)-\sum_{j=1}^m\hat{y}_j(t_i)\ell _j(x)\Big|f_X(x)\dx\nonumber \\
&=&\int_\R\Big|\sum_{j=1}^m\epsilon_j^A\ell _j(x)\Big|f_X(x)\dx.
\end{eqnarray}
Note that the $\ell_j(x)$ interpolation functions are identical as they depend solely on the $x$ values.
We arrive at the following error related to the ANNs,
\begin{eqnarray}
	\int_\R\Big|\sum_{j=1}^m\epsilon_j^A\ell _j(x)\Big|f_X(x)\dx
&\leq&\int_\R \sum_{j=1}^m\max\{|\epsilon_1^A|,\dots,|\epsilon_m^A|\}\ell _j(x)f_X(x)\dx\nonumber \\
&=&\int_\R \max\{|\epsilon_1^A|,\dots,|\epsilon_m^A|\}f_X(x)\dx\nonumber \\
&=&\max\{|\epsilon_1^A|,\dots,|\epsilon_m^A|\}
\end{eqnarray}
Considering the error introduced by SCMC in~\eqref{epsilon_scmc}, the total path wise error reads
\begin{eqnarray} \label{eq:pwerror}
    \E\left[ |g(x)-\tilde{g}_m(x) |\right] &\leq& \E\left[ |g(x)-g_m(x) |\right] + \E\left[ |g_m(x)-\tilde{g}_m(x) |\right] \nonumber \\ &\leq& \sqrt{|\epsilon_m|}+ \max\{|\epsilon_1^A|,\dots,|\epsilon_m^A|\}.
\end{eqnarray}
In other words, the expected pathwise error can be bounded by the approximation CDF error $\sqrt{|\epsilon_m|}$ plus the largest difference in the ANN approximated collocation points.

\subsubsection{Kolmogorov-Smirnov Test}
The Kolmogorov-Smirnov test, calculating the supremum of a set of distances, is used to measure  the nonparametric distance between two empirical cumulative distribution functions.
We perform the two-sample Kolmogorov-Smirnov test, as follows,
$$KS = \sup_z|F_Y(z) -\hat F_Y(z)|,$$
where $\hat F_Y(\cdot)$ and $F_Y(\cdot)$ are two empirical cumulative  distribution functions, one from the 7L-CDC solution and the other one from the reference distribution. We take the analytic solution of the GBM as the reference distribution.
\begin{remark}[Time horizon for 7L-CDC ]

	The information in Table~\ref{table:ann_train_data} is used to train the mapping function  between a realization (including marginal SC points) and its conditional SC points, via Equation~\eqref{eq:cijk}. For the marginal SC points in Equation~\eqref{yij}, however, we need training data up to terminal time $T$. So, we generate a second data set in which the time reaches $\tau_{max}$ (the terminal time of interest) and the upper value for $Y_0$ equals 5.  These two data sets are merged into one set in order to train the ANNs for the 7L-CDC methodology.
\end{remark}

\begin{figure}[!h]
\centering
\subfloat[KS statistic]{\label{fig:cdc_kstest_stat} {\includegraphics[width=0.5\textwidth]{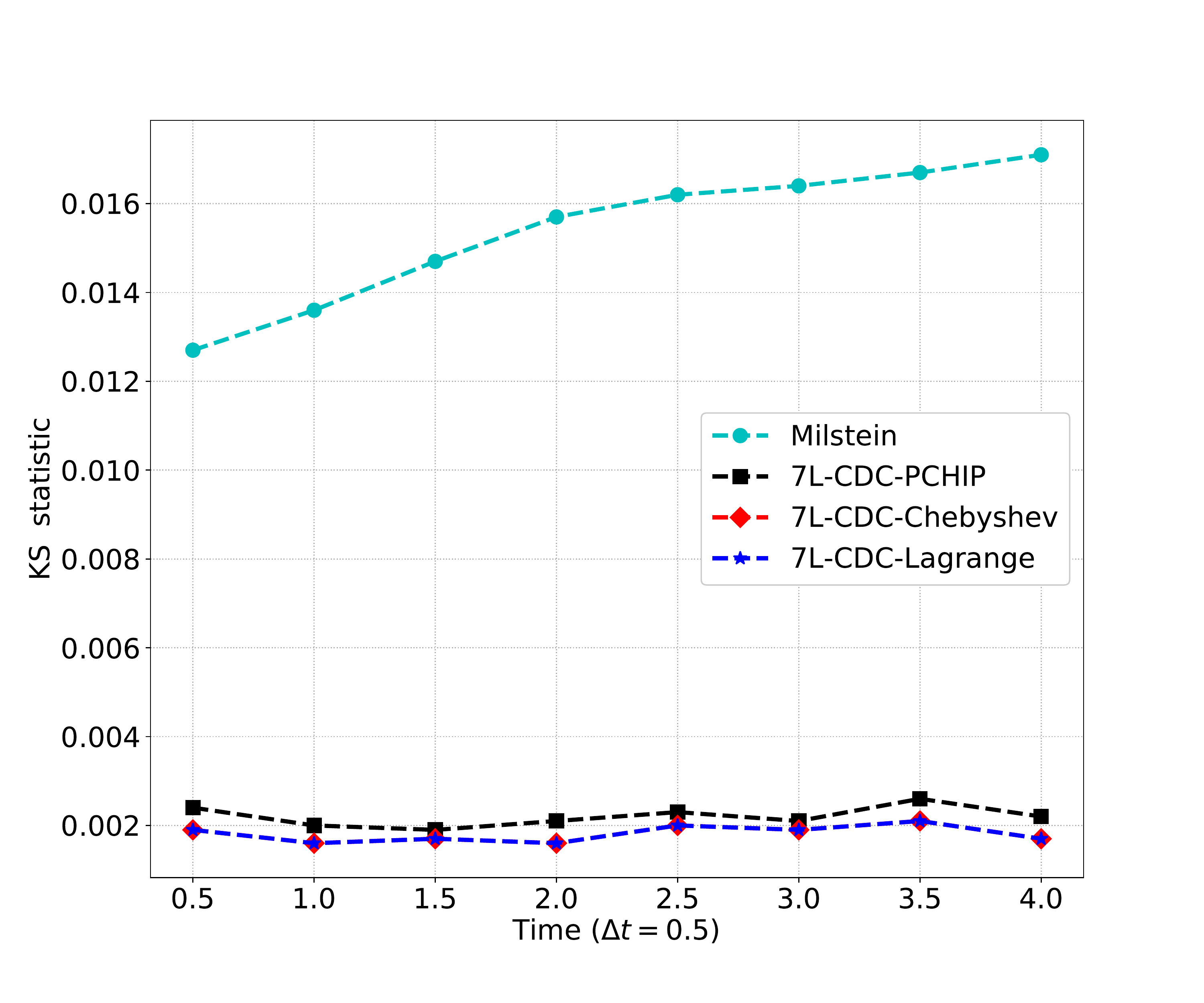}}}%
\subfloat[$P$-value]{\label{fig:cdc_kstest_pval}{\includegraphics[width=0.5\textwidth]{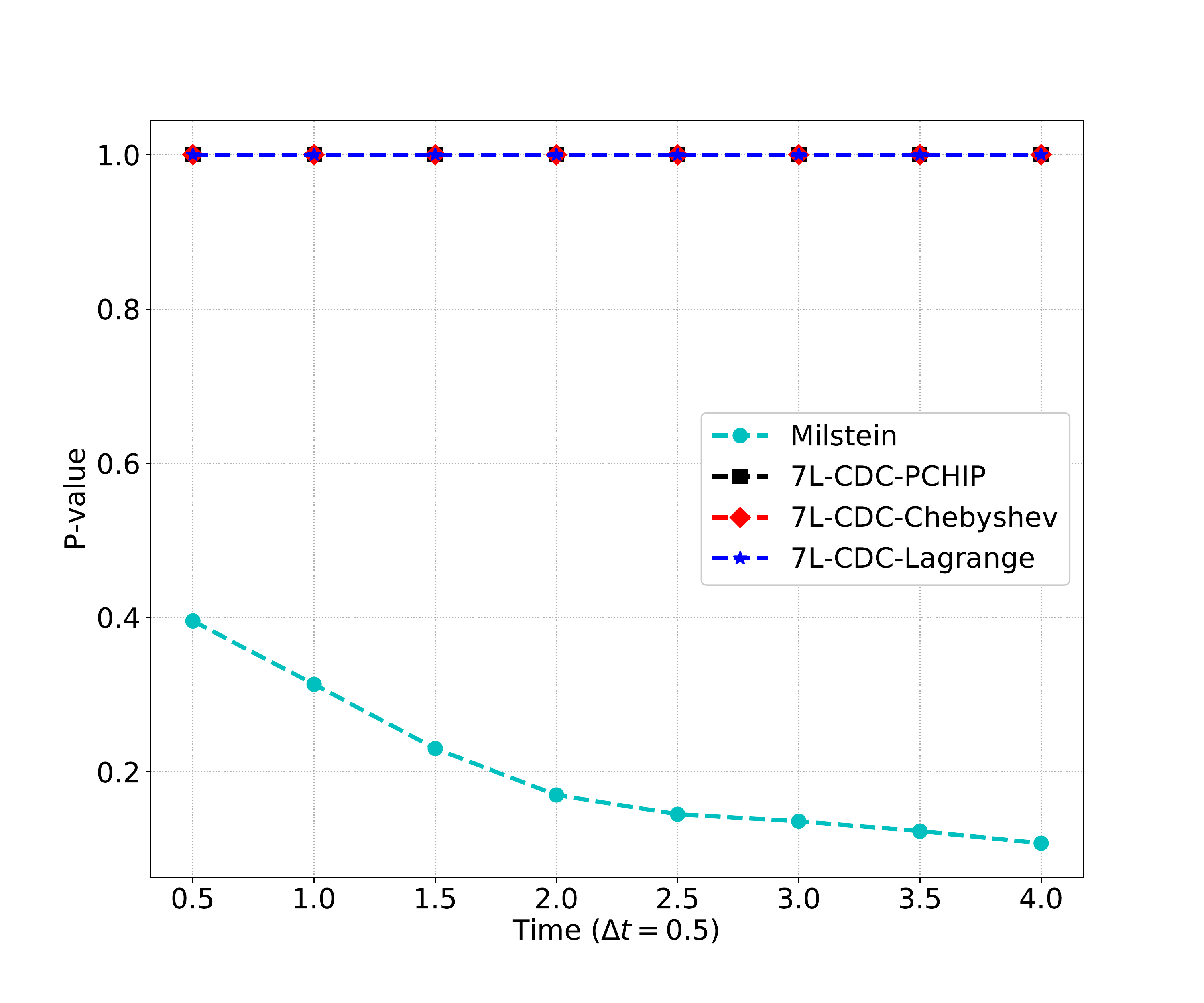}}}%
\caption{ The Kolmogorov-Smirnov test: $\Delta t =0.5, \mu=0.1, \sigma=0.3, Y_0=1.0$, with 10000 samples. When we have a small KS statistic or a large $P$-value,  the hypothesis that the distributions of the two sets of random samples are the same can not be rejected.}
\label{fig:ks-test}
\end{figure}

Figure~\ref{fig:ks-test} shows the Kolmogorov-Smirnov test at different time points based on 10000 samples.
We focus on the CDC methodology here, and compare the accuracy with the different interpolation methods in the figure. Clearly,  the KS statistic and also the corresponding $P$-values for the 7L-CDC schemes are much better than those of the Milstein scheme in Figure~\ref{fig:ks-test}.
This is an indication that the CDFs that originate from the 7L-CDC schemes resemble the target CDF much better, with high confidence. In addition, unlike the Milstein scheme
the 7L-CDC schemes exhibit an almost constant difference between the approximated and target CDFs with increasing time.

We will also analyze the costs of the different interpolation methods within 7L-CDC. The two steps which require interpolation are the computation of the conditional collocation points and the generation of conditional samples. The computational speed of the 7L-CDC scheme depends on the employed interpolation method, see Table~\ref{table:computing_time_interp}. In general, to generate a solution with the same strong order in the numerical error, the Milstein scheme will require more computation time,  here about 27 seconds,   while the 7L scheme needs 13 seconds and 7L-CDC (Barycentric version) 5 seconds when $\Delta t=1.0$. The larger the time step, the more computation time will be saved.

\begin{table}[!h]
\begin{center}
\caption{ The CPU running time (seconds) to reach the same accuracy (CPU: E3-1240, 3.40GHz): simulating 10,000 sample paths until terminal time $T=4.0$, based on $5\times5$ marginal/conditional SC points. Here, for the 7L scheme, PCHIP is used as the interpolant $g_m(\cdot)$ in Step 3 of Algorithm~I. }
\begin{tabular}{ c|c c c | c c c}
\hline
  \multirow{2}{*}{Method \slash Time (Sec.)} &  \multicolumn{3}{c|}{$\Delta t=1.0$} &  \multicolumn{3}{c}{$\Delta t=2.0$} \\
  & Create $C$ & Decom. $C$ & Total & Create $C$ & Decom. $C$ & Total
  \\ \hline
7L-CDC Barycentric & 0.054 & 4.93 & 4.98 & 0.027 & 2.48 & 2.51 \\
7L-CDC Chebyshev & 0.054 & 9.78 & 9.83 & 0.027 & 4.93 & 4.96 \\
7L-CDC PCHIP & 0.054 & 11.39 & 11.44 & 0.027 & 5.73 & 5.76 \\
\hline
7L scheme  & - & -& 12.80 & - & -& 6.39\\
Milstein & - & - & 27.01 & - & - & 27.70\\
\hline
\end{tabular}
\label{table:computing_time_interp}
\end{center}
\end{table}

Note that, in order to achieve a similar accuracy in the strong sense, the Euler-Maruyama  scheme requires a much finer time grid, by a factor of $\kappa=\Delta t/\Delta \tau$, than the 7L scheme. When $\kappa$ is sufficiently large, the 7L-CDC scheme outperforms the Euler-Maruyama scheme, in terms of both accuracy and speed. For example, in  Table~\ref{table:computing_time_interp}, $\kappa=100$ when $\Delta t =1.0$, and  $\kappa=200$ when $\Delta t =2.0$.  In other words, the ``on-line version'' of the Euler-Maruyama discretization is computationally slower than the on-line phase of the 7L scheme to achieve the same accuracy.  Additionally, computational time of the 7L scheme can be further reduced by
parallelization, for example, using GPUs .

\subsection{Path-wise Error Convergence} \label{sec:pathwise error}

In this section, we compare the path-wise errors of our proposed novel discretization with those of the classical discretization schemes.

\subsubsection{GBM process}
\begin{figure}[!h]
    \centering
    \includegraphics[width=0.6\textwidth]{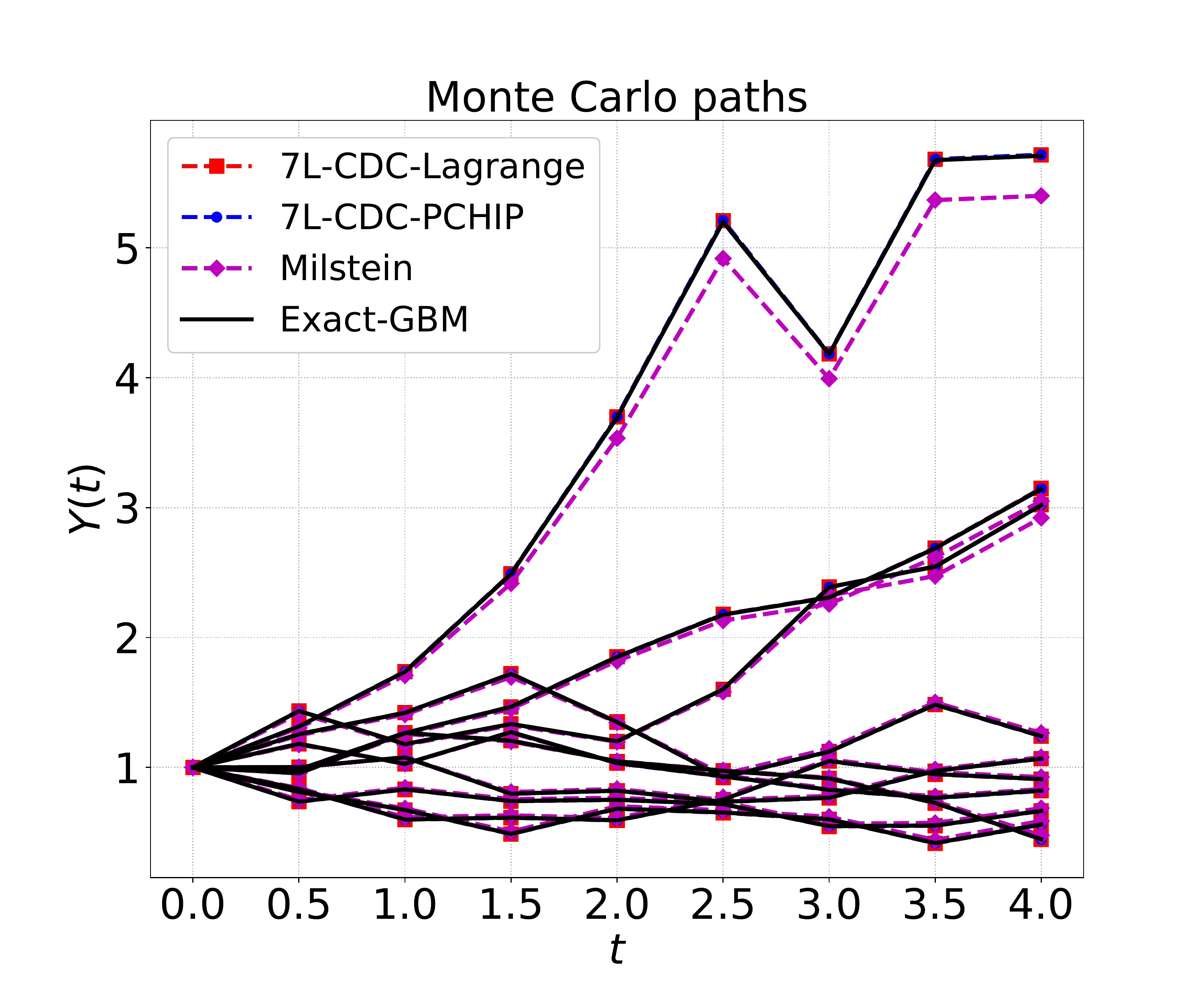}
    \caption{Paths generated by 7L-CDC: time step $\Delta t=0.5$, GBM with $\sigma=0.3$, $r=0.1$, $Y_0 =1.0$.  The paths are with Chebyshev interpolation, which are not plotted,  are identical to ones from Lagrange in this case. }
    \label{fig:paths_gbm_7Lcdc}
\end{figure}

We analyze here the strong convergence properties of the new methodology for the GBM process.
For GBM, the exact path is given by the expression~\eqref{exa}. The random number, which is drawn from $X\sim N(0,1)$, is the same for the exact solution ~\eqref{exa}, the novel schemes~\eqref{eq:condsam} and the Milstein scheme~\eqref{eq:milstein}.
The path-wise differences between the numerical schemes and the exact simulation are plotted in Figure~\ref{fig:paths_gbm_7Lcdc}.  When $\Delta t =0.5$, the 7L-CDC scheme presents superior paths as compared to the Milstein scheme, in terms of its path-wise error comparing to the exact path.

\begin{figure}[!h]
\centering
\subfloat[Strong convergence]{\label{fig:cdc_strong_converg}{\includegraphics[width=0.5\textwidth]{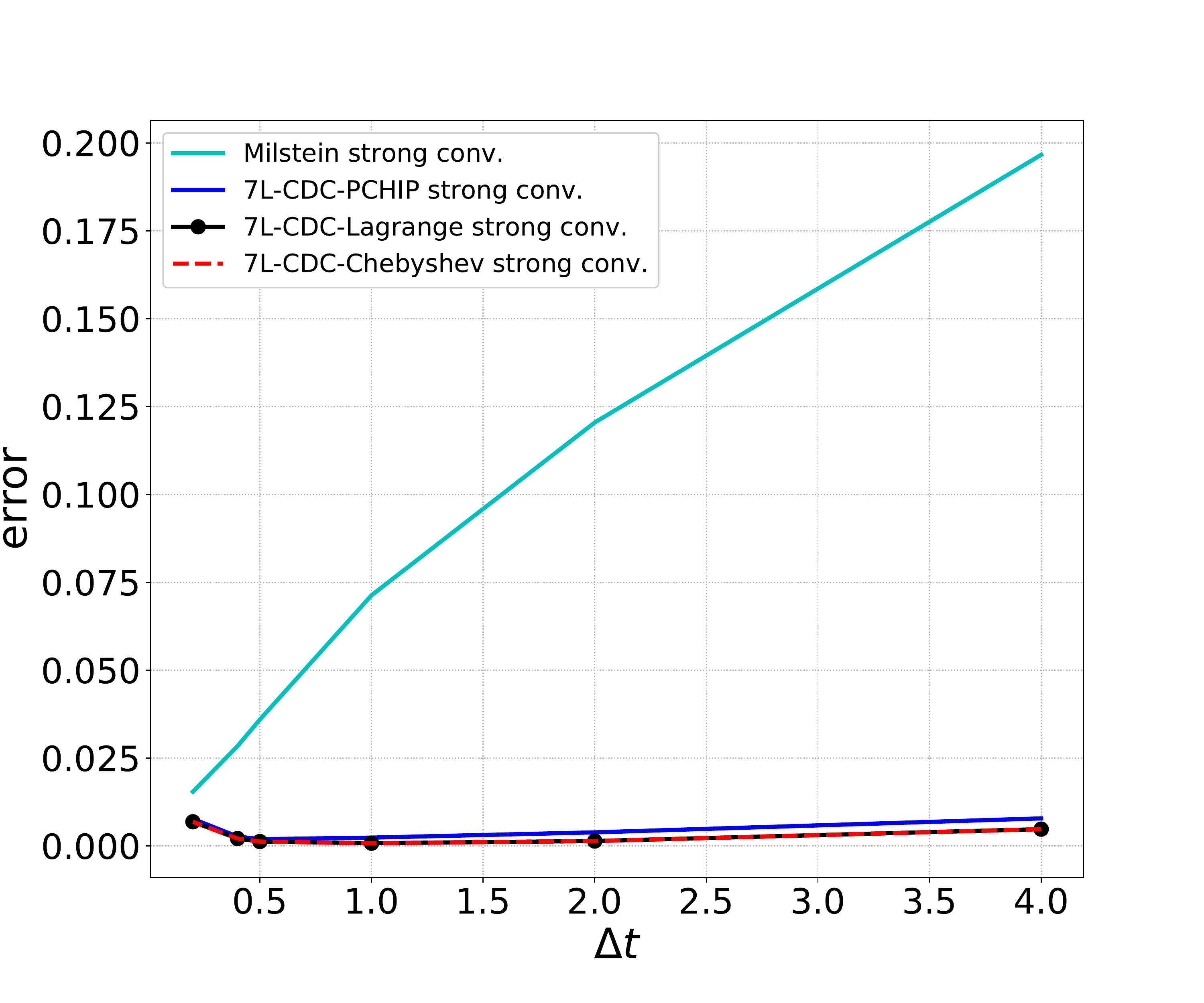}}}%
\subfloat[Weak convergence]{\label{fig:cdc_weak_converg} {\includegraphics[width=0.5\textwidth]{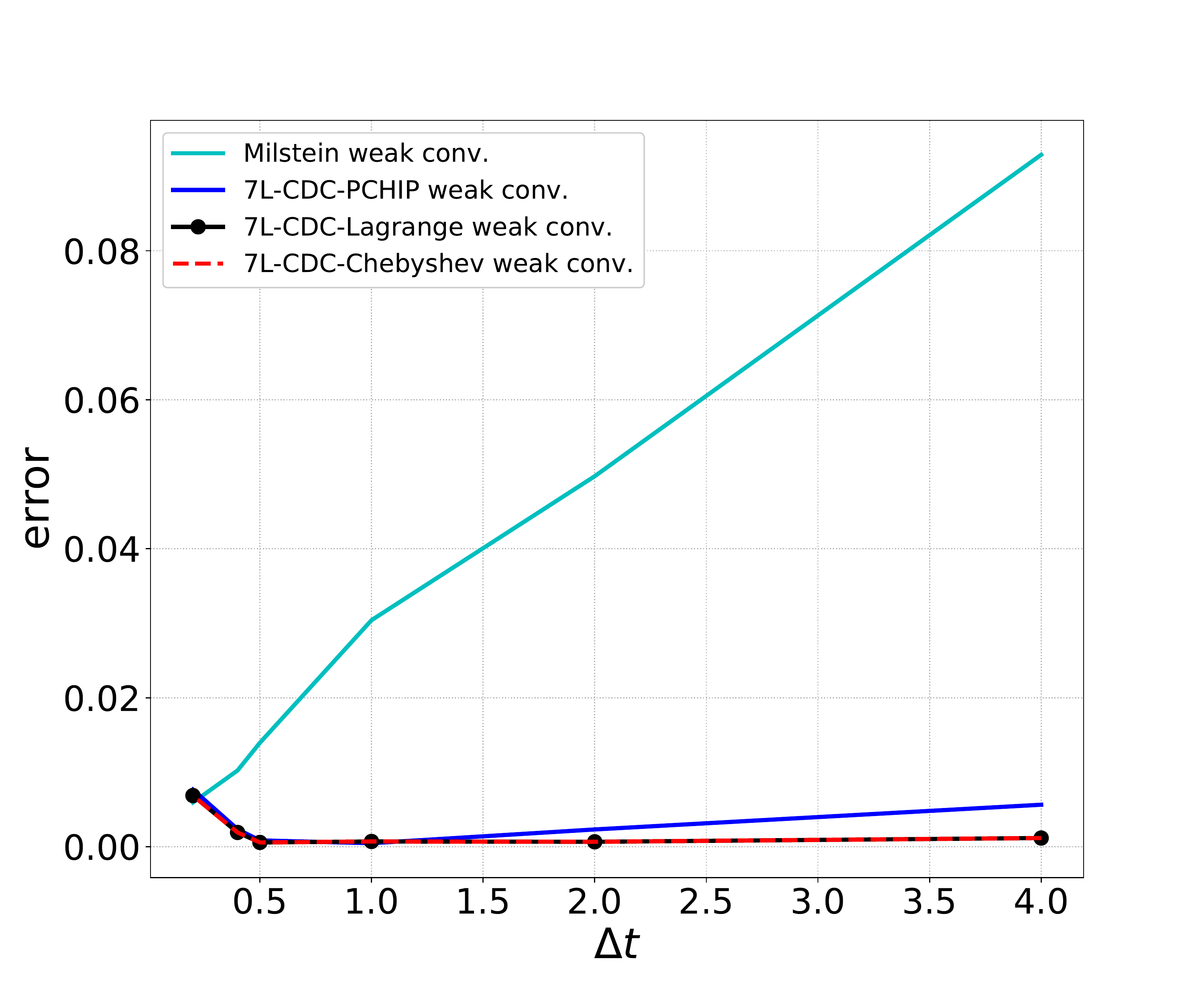}}}%
	\caption{The strong error is estimated as $\frac{1}{M}\sum |\check{Y}_k(T)-\hat{Y}_k(T)|$,  see Equation~\eqref{sstr} and the weak error by $\frac{1}{M}(\sum \check{Y}_k(T)-\sum\hat{Y}_k(T))$, see~\cite[][page 261]{finbook2019} for details on the computation of the convergence rate. There are $M=1000$ sample paths in total.}
\label{figure:convergence_error}
\end{figure}

As shown in Figure~\ref{figure:convergence_error},  the 7L-CDC scheme gives rise to flat, almost constant, strong and weak error convergence curves for many different $\Delta t$-values, suggesting a small, constant convergence error even with large time steps $\Delta t$. The Milstein scheme  has the strong order of convergence $O(\Delta t)$, so that a larger time step gives rise to a larger error.  When the time step becomes small,  more time points are needed to reach a time $T$, and then the resulting recursive error of the 7L-CDC scheme  increases.

The number of conditional collocation points, by which the conditional distribution at a next time point is mostly determined, has a significant contribution to the convergence order of the 7L-CDC scheme. As mentioned, we found  empirically that five conditional collocation points are preferable in terms of computing effort versus accuracy. CDC matrix $C$ is then of size $N\times5\times5$, that is, at each time point, there are five collocation points and each of these has five conditional collocation points.

\subsubsection{Ornstein-Uhlenbeck process}
\label{ouou}
Any SDE which can be solved by the Euler-Maruyama discretization can be solved by our ANN methodology, with improved strong convergence properties.
We also wish to confirm the strong convergence properties for another stochastic process in this section.

The mean reverting Ornstein-Uhlenbeck (OU) process~\citep{OUP1930} is defined as,
\begin{equation} \label{eq:oup}
\d Y(t) = -\lambda (Y(t)-\overline{Y}) \dt + \sigma \dW(t),  \;\;\;\;\; 0 \leq t \leq T,
\end{equation}
with $\overline{Y}$ the long term mean of $Y(t)$, $\lambda$ the speed of mean reversion, and $\sigma$ the volatility. The initial value is $Y_0$, and the model parameters are $\boldsymbol\theta :=\{\overline{Y},\sigma,\lambda\}$.   Its analytical solution is given by,
\begin{equation} \label{eq:oup-solu}
    Y(t) \stackrel{d}{=} Y_0 e^{-\lambda t} + \overline{Y}(1-e^{-\lambda t}) +\sigma \sqrt{ \frac{1-e^{-2\lambda  t}}{2\lambda}} X,
\end{equation}
with $t_0=0$, $X \sim \mathcal{N}(0,1)$.  Equation~\eqref{eq:oup-solu} is used to compute the reference value to the path-wise error and the strong convergence.

We employ the same data-driven procedure as for GBM  to discretize and solve the OU process.   In the training phase, the Euler-Maruyama scheme~\eqref{eq:euler} is used to discretize the OU dynamics and generate the data set. Note that the Milstein and Euler schemes are identical in the case of the OU process. As the OU process is a Markov process,  we again can vary $Y_0$ to find the relation between the conditional SC points and the marginal SC points (i.e. as in Equation~\eqref{eq:cijk}).  Similar to Table~\ref{table:ann_train_data}, we employ five SC points to learn within the ANN,
with $\Delta \tau =0.01$, $\tau_{max}=4.1$, $N_\tau=500$, $M_L=410$, see Section~\ref{sec:trainprocess} for the details of the training process.

\begin{figure}[!h]
\centering
\subfloat[Path-wise error ($\Delta t=1.0$)]{\label{fig:paths_OUP_cdc} {\includegraphics[width=0.5\textwidth]{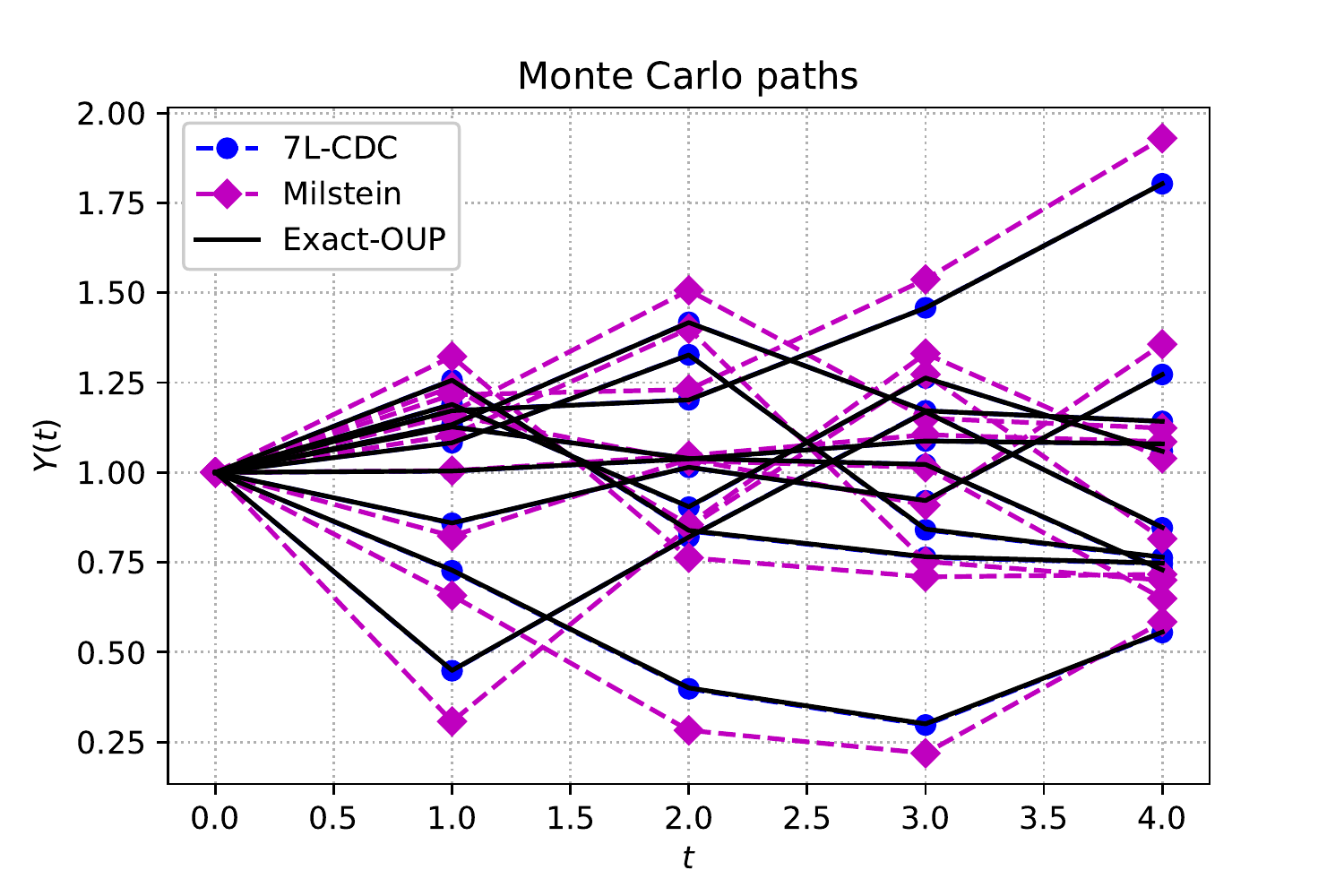}}}
\subfloat[Strong convergence]{\label{fig:anncdc_strong_converg_oup}{\includegraphics[width=0.5\textwidth]{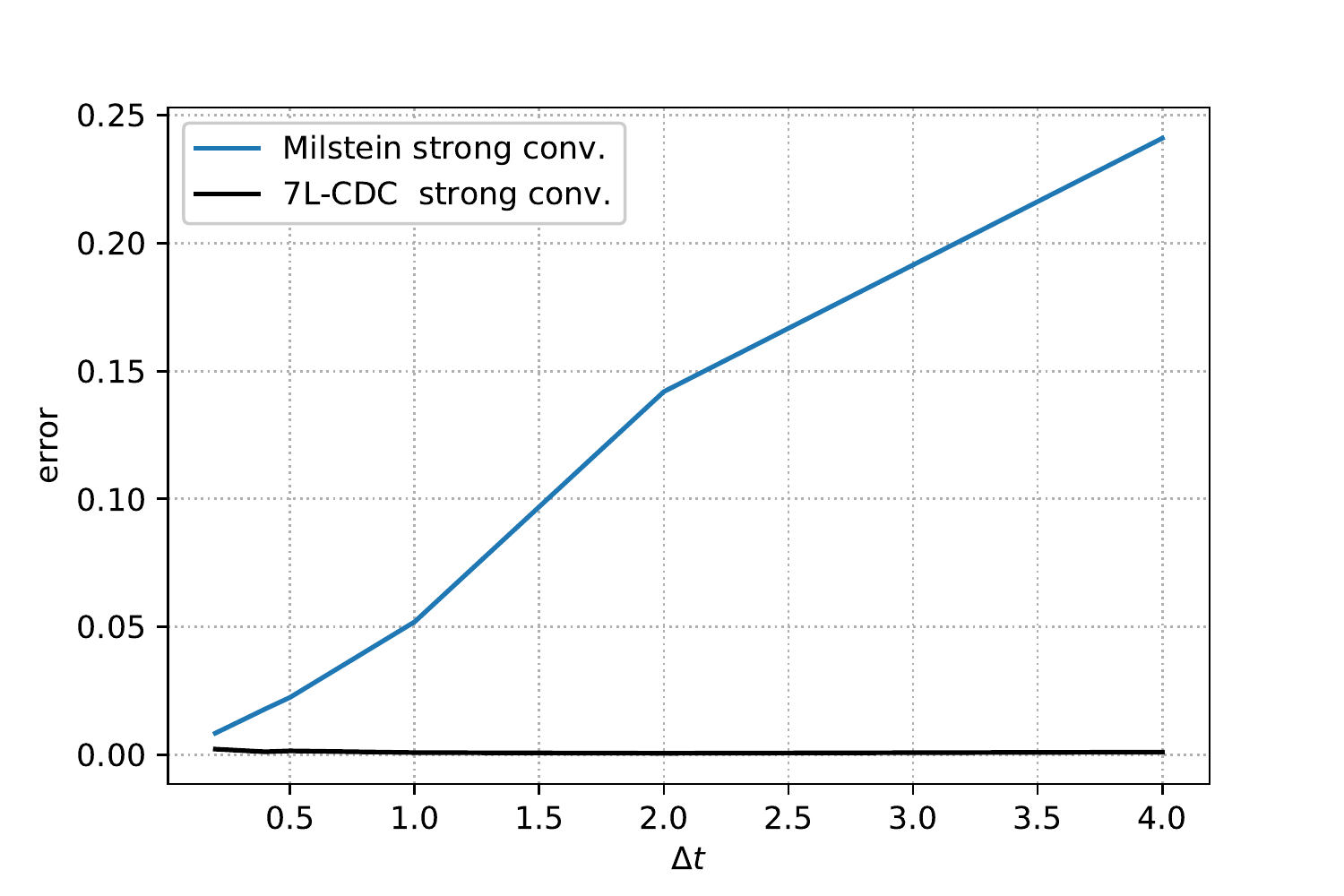}}}
\caption{Paths and strong convergence for the  OU process, using $\lambda=0.5$, $ \overline{Y}=1.0$, $\sigma=0.3$,  $Y_0=1.0$. The sample paths with barycentric, Chebyshev and PCHIP interpolation overlap for the 7L-CDC scheme.  There are five marginal and five conditional SC points at each time point.}
\label{fig:oup-test}
\end{figure}

After the training, the obtained ANNs will be applied to solve the OU process with specific parameters and details of our interest. We provide an example in Figure~\ref{fig:oup-test}, which confirms that the sample paths generated by 7L-CDC are as accurate as the exact solution, and the error, in the sense of strong convergence, stays close to zero even with a large time step.

\subsection{Applications in Finance}
The possibility to take large time steps and still get accurate SDE solutions, is certainly interesting in computational finance, as there are several financial products that are updated on a daily basis (think of an over-night interest rate), whereas monitoring of financial contracts and risk management monitoring is typically only done on a weekly, monthly of even yearly basis. In such situations, our novel scheme will be useful. Research into large time step simulations is state-of-the-art in computational finance, see the exact (and almost exact) Monte Carlo simulation papers, like~\cite{KayaBroadie2006,LEITAO2017461} for the SABR and Heston stochastic volatility asset dynamics, respectively.

\subsubsection{The Asian option}
Moreover, the strong convergence property of an SDE discretization is important in many cases.  When valuing so-called path-dependent options, for example, improved strong convergence enhances the convergence of a Monte Carlo simulation. Options are governed by their pay-off function (i.e. the option value at the final time of the contract, $t=T$).
Here we consider a path-dependent exotic option,  the so-called European-style Asian option, which has a payoff that is based on a time-averaged underlying stock price. For example, the pay-off of a fixed strike Asian option is given by
$$
V_A(T) = \max{(A(T)-\tilde{K}, 0)},
$$
where  $T$ is the option contract's expiry time, and $\tilde{K}$ is the predetermined strike price. Here $A(T)$ denotes the discrete arithmetic average of the stock prices over $N_b$ monitoring dates $\{t_1,\ldots, t_k\} \in [0,T]$,
$$A(T) = \frac{1}{N_b} \sum_{k=1}^{N_b} \hat{Y}(t_k), $$
where $\hat{Y}(t_k)$ is the observed stock price at time $t_k$, $ 1 \leq t_k \leq T$.
Averaging thus takes place in the time-wise direction, and we consider pricing financial options based on the discrete arithmetic average of a number of stock prices.

We assume here that the underlying stock price follows Geometric Brownian motion, as in Equation~\eqref{eq:gbmq}, under the risk-neutral measure, meaning $\mu\equiv r$, where $r$ is the risk-free interest rate. There is a cash account $M(t)$, governed by $\d M(t) = r M(t)\dt$.
The value of  European-style Asian option is then given by
\begin{equation}
	V_A(t) = \e^{- r (T-t)} \mathbb{E}^{\mathbb{Q}} \biggl[ \max ( A(T)-\tilde{K} , 0) \bigg|  \F(t)\biggr].
\end{equation}
Because the pay-off is clearly a path-dependent quantity for such options, it is expected that an improved strong convergence, obtained with the variant 7L-CDC, will result in superior convergence, as compared to classical numerical discretization schemes.
\begin{table}[!h]
\begin{center}
  \caption{Pricing Asian European-style option with a fixed strike price, using $Y_0=1.0$,  $\tilde{K}=Y_0$, $r=0.1$, $T=\Delta t \times N_b$.}
   \begin{tabular}{  c | c | c  | c }
    \hline
    & method &  $\Delta t =1.0$, $N_b$=4 & $\Delta t=0.5$, $N_b$=8  \\
   \hline
   \multirow{3}{*}{$\sigma$=0.30}
    & Analytic MC  & 0.24886257 (0.00\%)  & 0.22403982 (0.00\%)\\
    & Milstein MC & 0.23077000 (7.27\%)   & 0.21558276 (3.77\%) \\
    & 7L-CDC & 0.24871446 (0.06\%)   & 0.22404571 (0.00\%)\\
    \hline
    \multirow{3}{*}{$\sigma$=0.40}
    & Analytic MC  & 0.28515109 (0.00\%)  & 0.25723594 (0.00\%) \\
    & Milstein MC  & 0.26394277 (7.44\%)  & 0.24717425 (3.91\%)  \\
    & 7L-CDC & 0.28482371 (0.11\%)  & 0.25647592 (0.30\%)\\
    \hline
  \end{tabular}
\label{table:perfomance_asian_pirce}
\end{center}
\end{table}

The relative error is presented, which is defined as $$\epsilon_{rel} = \biggl|\frac{V_A^{ref}(t_0)-V_A(t_0)}{V^{ref}_A(t_0)}\biggr|,$$  where $V^{ref}_A(t_0)$ is based on the exact GBM Monte Carlo simulation. As shown in Table~\ref{table:perfomance_asian_pirce}, the 7L-CDC scheme gives highly accurate Asian option prices, compared to the Milstein scheme. As the accuracy of Asian option prices depends directly on the accuracy of the realized paths, an increasing number of monitoring dates will give rise to higher
accuracy by 7L-CDC.

Next, we focus on the Asian option's sensitivity.
The sensitivity of the option price with respect to volatility $\sigma$ is called vega, which can be computed in a pathwise fashion (see Chapter 7 in~\cite{glasserman2004monte}), as follows,
\begin{equation} \label{eq:Vwrtvol}
    \frac{\partial V}{ \partial \sigma}=e^{-rT}\mathbb{E}^{\mathbb{Q}} \left[ \sum_{i=1}^N \frac{\partial V(T,Y(t_i);\sigma)}{\partial Y(t_i)} \frac{\partial Y(t_i)}{\partial \sigma} \big| Y_0  \right].
\end{equation}
The  chain rule is employed to derive the sensitivity. First of all, we compute the gradient of the payoff function with respect to the underlying stock price, by
\begin{equation}
    \frac{\partial V(T,Y(t_i))}{ \partial Y(t_i)} = \frac{1}{N} \mathbbm{1}_{A(T)>\tilde{K}}.
\end{equation}
Then,  the derivative of the stock price at time $t_i$ with respect to the model parameter, $\frac{\partial Y(t_i)}{\partial \sigma}$, can be found with the trained ANNs, as given by Equation~\eqref{eq:ywrttheta}. Vega can be estimated by,
\begin{equation} \label{eq:Vwrtvolann}
    \frac{\partial V}{ \partial \sigma} \approx e^{-rT}\frac{1}{N}\mathbb{E}^{\mathbb{Q}} \left[ \sum_{i=1}^N \left( \mathbbm{1}_{A(T)>\tilde{K}}  \sum_{j=0}^{m-1} \frac{\partial \hat{H}_j}{\partial \sigma} p_j(X) \right) \big| Y_0 \right].
\end{equation}

When there are $M$ sample paths, we have,
\begin{equation} \label{eq:Vwrtvolpaths}
    \frac{\partial V}{ \partial \sigma} \approx e^{-rT} \frac{1}{M} \frac{1}{N} \left[ \sum_{k=1}^{M} \sum_{i=1}^N \left( \mathbbm{1}_{A(T)>\tilde{K}}   \sum_{j=0}^{m-1} \frac{\partial \hat{H}_j}{\partial \sigma} p_j(\hat{X}_{k,i+1}) \right) \big| Y_0  \right].
\end{equation}

\begin{figure}[!h]
\centering
\subfloat[Path-wise sensitivity]{\label{fig:paths_sense_gbm_scmc}{\includegraphics[width=0.5\textwidth]{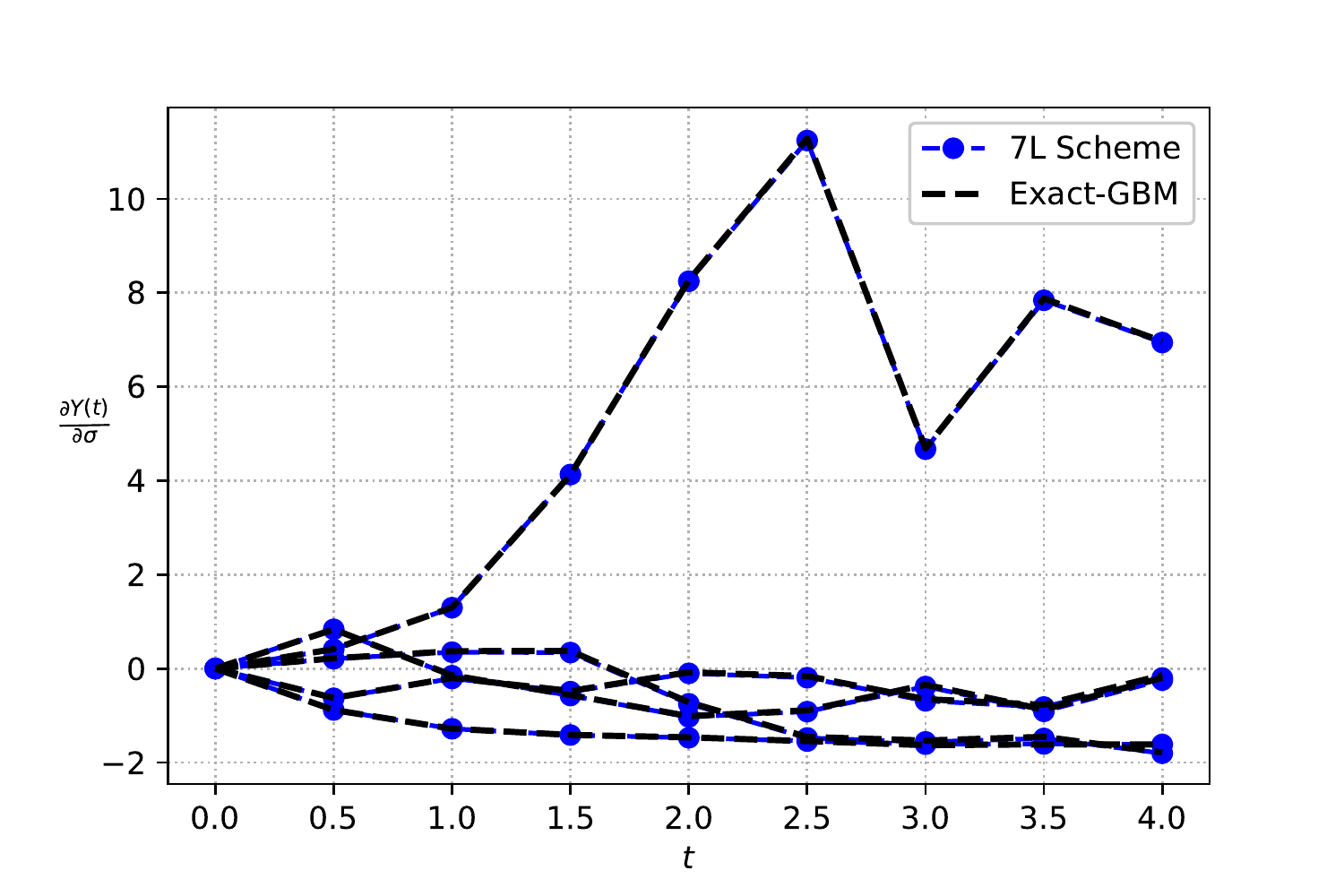}}}
\subfloat[Path-wise vega]{\label{fig:pathwise-vega} {\includegraphics[width=0.5\textwidth]{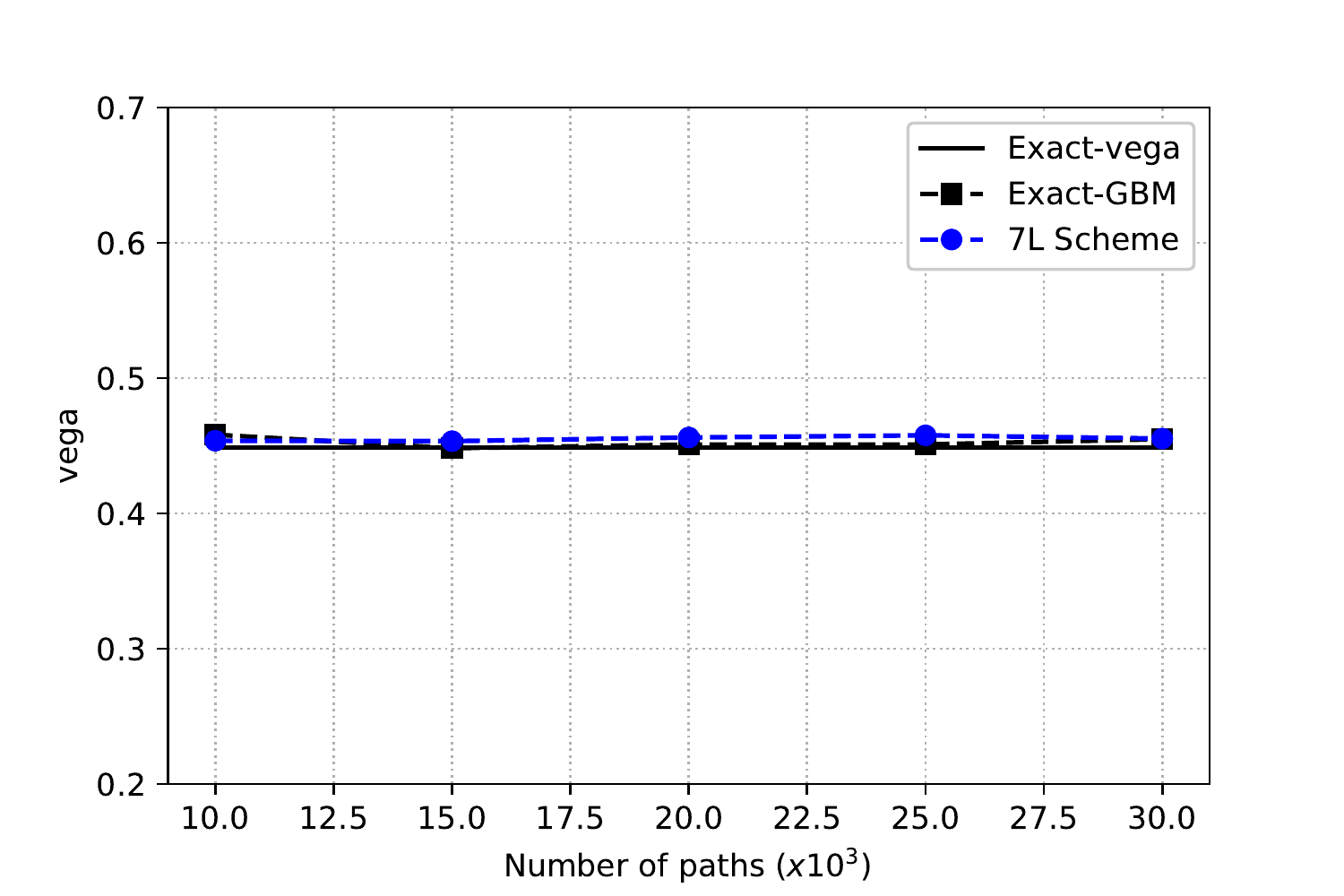}}}
	\caption{ Path-wise estimator of vega: Exact vega is calculated by means of the central finite difference. The parameters are $Y_0=1.0$, $r=0.05$, $\tilde{K}=Y_0$,  $\sigma=0.3$, $\Delta t =1.0$, $N_b=4$, $T=\Delta t \times N_b=4.0$.}
\label{figure:Path-wise-sensitivy}
\end{figure}

As the realization $\hat{Y}_i$ is a function of the model parameters,
at time $t_{i+1}$,  the derivative with respect to the volatility in Equation~\eqref{eq:ywrtthetapath} becomes
\begin{equation} \label{eq:gbmscmcpathwise}
    \frac{\partial \hat{H}_{j}(\hat{Y_i},\sigma)}{\partial \sigma} = \frac{\partial \hat{H}_{j}(\sigma;\hat{Y_i})}{\partial \sigma} + \frac{\partial \hat{H}_{j}(\hat{Y_i};\sigma)}{\partial \hat{Y_i}}   \frac{\partial \hat{Y}_i}{ \partial \sigma},
\end{equation}
where $\frac{\partial \hat{Y}_i}{ \partial \sigma}$ is known at the previous time point. Like simulating the Monte Carlo paths, the calculation of this derivative is done  iteratively.
Figure~\ref{fig:paths_sense_gbm_scmc} compares the path-wise sensitivities obtained via Equations~\eqref{eq:gbmsense} and~\eqref{eq:gbmscmcpathwise}. Clearly, the path-wise derivative by the 7L scheme is very similar to the analytical solution.  Figure~\ref{fig:pathwise-vega} confirms that the ANN methodology computes a highly accurate Asian option vega by means of the above path-wise sensitivity. Summarizing, the sensitivity with respect to model parameters can highly accurately be obtained from the trained ANNs. As the 7L-CDC scheme is composed of marginal and conditional collocation points, the above procedure of computing the path-wise sensitivity is also applicable to the variant 7L-CDC, by using the chain rule.

\subsubsection{Bermudan option valuation}  \label{sec:bermudan}
When dealing with so-called Bermudan options, the option contract holder has the right (but not the obligation) to exercise the option contract at a finite number of pre-specified dates up to final time $T$.
At an exercise date, when the holder decides to exercise the Bermudan option, she immediately obtains the current payoff value of the contract. Alternatively, she may also wait until the next exercise opportunity.
The Bermudan option can be exercised at the following set of exercise dates, $\{t_0,t_1,...,t_{N_b}\}$, with a constant time difference, $\Delta t = t_{i}-t_{i-1}$, for any $0<i\leq N_b$.

In this experiment, we compare the performance of the new 7L-CDC discretization scheme with a classical scheme. Valuation of the Bermudan option will take place by means of the well-known Longstaff-Schwartz Monte Carlo method (LSMC)~\citep{Longstaff}, a least squares Monte Carlo method. 
The Longstaff-Schwartz algorithm is presented, for convenience,  in the appendix.

The difference between a large time step simulation and a classical simulation, like the Milstein scheme, is that a classical scheme requires additional time steps to be taken between the early-exercise dates of the Bermudan option, while with the
7L-CDC scheme, we can perform  one-step Monte Carlo simulation without any intermediate grid points between adjacent early-exercise dates.

We also assume here that the underlying stock price follows Geometric Brownian motion, as in Equation~\eqref{eq:gbmq}, under the risk-neutral measure, with $\mu\equiv r$. 
A Bermudan put option, with risk-free interest rate $r=0.1$,  pay-off function $V(t_j) =  \max{(\tilde{K}-Y(t_j),0)}$ with strike price $\tilde{K}=1.1$ and initial stock price $Y_0=1.0$, is priced based on $M=100,000$ Monte Carlo paths.  The matrix size within the 7L-CDC scheme is set to $N_b\times 5\times5$.
The terminal time is $T=\Delta t \times M_B$ with a constant time step $\Delta t$. The random seed is chosen to be zero when drawing random numbers. We compare the relative errors $|\frac{V^{ref}(t_0)-V(t_0)}{V^{ref}(t_0)}|$,  where $V^{ref}(t_0)$ is computed with the help of a Monte Carlo method based on the exact simulation of GBM~(\ref{exa}).

\begin{table}[!h]
\begin{center}
  \caption{Bermudan put option prices based on large time step Monte Carlo simulations.}
   \begin{tabular}{  c | c | c | c | c }
    \hline
    & method &  $\Delta t =1.0$, $N_b$=4 & $\Delta t =0.5$, $N_b$=4 & $\Delta t=0.5$, $N_b$=8  \\
   \hline
   \multirow{3}{*}{$\sigma$=0.30}
    & Analytic MC  & 0.15213858(0.00\%) & 0.14620214(0.00\%) & 0.16161876(0.00\%)\\
    & Milstein MC & 0.13872771(8.81\%)  & 0.14065252(3.80\%)  & 0.15429369(4.53\%) \\
    & 7L-CDC & 0.15234901(0.14\%)  & 0.14648443(0.19\%)  & 0.16196264(0.21\%)\\
    \hline
    \multirow{3}{*}{$\sigma$=0.40}
    & Analytic MC  & 0.21459038(0.00\%) & 0.19552454(0.00\%) & 0.22340304(0.00\%) \\
    & Milstein MC  & 0.19598488(8.67\%) & 0.18790933(3.89\%) & 0.21297732(4.67\%)  \\
    & 7L-CDC & 0.21474619(0.07\%) & 0.19590733(0.20\%)  & 0.22389360(0.22\%)\\
    \hline
  \end{tabular}
\label{table:perfomance_bermudan_pirce}
\end{center}
\end{table}

As shown in Table~\ref{table:perfomance_bermudan_pirce}, the option prices based on the 7L-CDC Monte Carlo simulation are highly satisfactory, and the related error does not increase with larger time steps $\Delta t$. In contrast, a larger time step gives rise to significant pricing errors, in the case of the Milstein discretization.

\begin{remark}
In principle, a sample value $\hat{Y}_i$ can be any rational number. So, a path value may reach a larger stock price than the prescribed upper bound in Table~\ref{table:ann_train_data}.
The stock prices outside the training interval are called {\em outliers}. Outliers did not appear in the experiments of  Table~\ref{table:perfomance_bermudan_pirce}.  As an alternative method to avoid the appearance of outliers, one may scale the asset price, to remove the dependence on the initial value. 	For example,  GBM can be scaled by $ \bar{Y}(t)= \frac{Y(t)}{Y(t_0) \e^{rt}}.$
Using It\^o's lemma, we have a drift-less process,
$ \d \bar{Y}(t) = \bar{Y}(t) \sigma \dW,$
where the initial value $\bar{Y}_0 =1.0$.
The following formula returns the original variable,
$    Y(t)= \bar{Y}(t)  Y(t_0)  \e^{rt}.$
In such case, scaling guarantees a fixed initial value, for example, $Y_0=1.0$.
\end{remark}

\section{Conclusions and Outlook}
\label{cenou}

We develop a data-driven numerical solver for stochastic differential equations,  by which large time step simulations can be carried out accurately in the sense of strong convergence. With a combination of  artificial neural networks and the stochastic collocation Monte Carlo method, a small number of stochastic collocation points are learned by the ANN to approximate a nonlinear function which can be used to compute the unknown collocation points. Theoretical analysis indicates that the numerical error is controllable and does not increase when the simulation time step increases.

There are several advantages to the proposed approach. The powerful expressive ability of neural networks enables the ANNs to accurately approximate stochastic collocation points. The compression-decompression method reduces the computational costs, so that the numerical method can be applied in practice. In finance, the proposed big time step methodology will be highly beneficial for the generation for path-dependent financial option contracts or in risk management applications.

As an outlook,  it will be relevant to extend the introduced methodology to solving higher-dimensional or more involved SDE dynamics. We will define multi-dimensional stochastic collocation points (i.e. by means of a tensor) for a multi-dimensional system of SDEs, and choose a Convolutional Neural Network~\citep{2015nature} to  efficiently process these collocation points. Of course, this may not trivially generalize to truly high-dimensional systems, but approximation of moderate dimensionality should be possible The computational speed can be further improved by parallel computation, for example, on GPUs. Non-Markovian processes may also be solved with a large time step by the proposed ANN method, where the conditional collocation points are dependent on past realizations.  Fractional Brownian motion~\citep{FBM1968} forms a relevant example, which is used for the simulation of rough volatility in finance~\citep{roughvol2018}.   In such a context, advanced variants of fully connected neural networks, e.g., recurrent neural networks (RNN) or long short-term memory (LSTM) networks (see a review in~\cite{RNNreview}), are recommended when approximating the nonlinear transition probability function, for example, Equation~\eqref{eq:condY}. 

As another outlook, Multilevel Monte Carlo (MLMC) methods, as developed by~\cite{mlmc2008,giles_2015}, form another interesting research topic for our large time step accurate discretisation schemes.
It is well-known that the strong convergence properties of SDE discretizations impact the efficiency of the MLMC methods.


\section{Acknowledgments}
S. Liu would like to thank the China Scholarship Council (CSC) for the financial support.

\bibliographystyle{plain}
\bibliography{main}

\appendix
\section{Longstaff-Schwartz Algorithm}
For convenience, we detail the Longstaff-Schwartz LSMC algorithm here.

\zbox{
	{\bf Algorithm: 7L scheme Longstaff-Schwartz algorithm} \\
	\label{alg3}
\begin{enumerate}
\item Divide the time horizon into $N_b$ intervals.
\item Simulate $M$ stock price paths $\hat{Y}_{i,j}$ ($0 \leq i \leq N_b$, $1 \leq j \leq M$), using the 7L-CDC methodology;

\item Price the Bermudan option by means of the Longstaff-Schwartz Monte Carlo method:

\begin{enumerate}
\item At terminal time $T_{N_b}$, calculate the payoff $\hat{V}_{N_b,j} = V\Big(Y_{N_b,j}\Big)$ for all paths $j$, where $V(\cdot)$ is the payoff function.
\item  Perform a backward recursion, from $i=N_b-1$ until $i=0$ as follows:

\item Compute the discounted continuation value at time $t_i$, i.e.,
\begin{equation} \label{eq:cvalue}
    \hat{\eta}_{i,j}:= \e^{-r\Delta t} \hat{V}_{i+1,j}
\end{equation}

\item Perform least squares regression at time $t_{i}$, based on the cross-sectional information $\hat{Y}_{i,j}$ and $\hat{\eta}_{i,j}$ to estimate the conditional expectation function,
\begin{equation} \label{eq:lsregression}
    \bar{\eta_i}(\hat{Y}) = \sum_{k=1}^{M_k} \beta_k B_k(\hat{Y})
\end{equation}
where $M_k$ is the number of basis functions $B_k(S)$ (polynomial basis, here, $M_k=3$), and the coefficients $\beta_k$ are constant over different paths $j$. Note that only in-the-money paths are considered in Equations~\eqref{eq:cvalue} and~\eqref{eq:lsregression},
\item For each path $j$, compare the immediate exercise value $V\Big(\hat{Y}_{i,j}\Big)$ with the estimated continuation value
$\bar{\eta_i}\Big(\hat{Y}_{i,j}\Big)$: If $V\Big(\hat{Y}_{i,j}\Big) \geq \bar{\eta_i}\Big(\hat{Y}_{i,j}\Big)$, then $\hat{V}_{i,j} = V\Big(\hat{Y}_{i,j}\Big)$; else  $\hat{V}_{i,j} = \bar{\eta_i}\Big(\hat{Y}_{i,j}\Big)$.
\end{enumerate}
\item Calculate the option price $V(t_0)$ at the initial time,
\begin{equation}
	V(t_0) = \frac{1}{M} \sum_{j=1}^M \hat{V}_{0,j}.
\label{eq:Price-LS}
\end{equation}
\end{enumerate}
}
\end{document}